%% file: 00-mssp.tex
\documentclass[preprint,3p,a4paper]{elsarticle}

\usepackage{setspace}
\usepackage{times}
\usepackage{pdflscape}
\usepackage{soul}
\usepackage{url}
\usepackage[hidelinks]{hyperref}
\usepackage[utf8]{inputenc}
\usepackage{graphicx}
\usepackage{amsmath}
\usepackage{amsthm}
\usepackage{mathtools}
\usepackage[ruled,vlined]{algorithm2e}
\usepackage{xcolor}
\usepackage{amssymb}
\usepackage{bm}
\usepackage{tikz}
\usepackage{booktabs}
\usepackage{multirow}
\usepackage{pgfplots}
\usepgfplotslibrary{groupplots}
\definecolor{marine}{RGB}{0,32,96}
\definecolor{navy}{RGB}{0,0,128}
\definecolor{maroon}{RGB}{128,0,0}
\definecolor{olivegreen}{RGB}{85,107,47}
\definecolor{gray}{RGB}{102,102,102}
\definecolor{green}{RGB}{5,66,8}
\definecolor{skyblue}{rgb}{0.3010, 0.7450, 0.9330}
\definecolor{purple}{rgb}{0.4940, 0.1840, 0.5560}
\definecolor{orange}{rgb}{0.9290, 0.6940, 0.1250}
\definecolor{brown}{RGB}{161,121,124}
\definecolor{deepblue}{rgb}{0.0, 0.0, 1.0}

\usepackage[final,authormarkup=none,commentmarkup=footnote]{changes}

\definechangesauthor[color=maroon]{AE}
\definechangesauthor[color=red]{YS}
\definechangesauthor[color=blue]{yunzhuang}
\definechangesauthor[color=orange]{XL}
\definechangesauthor[color=green]{ACE}

\usepackage[caption=false,font=normalsize,labelfont=sf,textfont=sf]{subfig}

\pgfplotsset{every tick label/.append style={font=\tiny}}
\pgfplotsset{
	box plot width/.initial=4em,
	box plot/.style={
		/pgfplots/.cd,
		black,
		only marks,
		mark=-,
		mark size=\pgfkeysvalueof{/pgfplots/box plot width},
		/pgfplots/error bars/.cd,
		y dir=plus,
		y explicit,
	},
	box plot box/.style={
		/pgfplots/error bars/draw error bar/.code 2 args={%
			\draw [line width=0.20mm]  ##1 -- ++(\pgfkeysvalueof{/pgfplots/box plot width},0pt) |- ##2 -- ++(-\pgfkeysvalueof{/pgfplots/box plot width},0pt) |- ##1 -- cycle;
		},
		/pgfplots/table/.cd,
		y index=2,
		y error expr={\thisrowno{3}-\thisrowno{2}},
		/pgfplots/box plot
	},
	box plot top whisker/.style={
		/pgfplots/error bars/draw error bar/.code 2 args={%
			\pgfkeysgetvalue{/pgfplots/error bars/error mark}%
			{\pgfplotserrorbarsmark}%
			\pgfkeysgetvalue{/pgfplots/error bars/error mark options}%
			{\pgfplotserrorbarsmarkopts}%
			\path ##1 -- ##2;
		},
		/pgfplots/table/.cd,
		y index=4,
		y error expr={\thisrowno{2}-\thisrowno{4}},
		/pgfplots/box plot
	},
	box plot bottom whisker/.style={
		/pgfplots/error bars/draw error bar/.code 2 args={%
			\pgfkeysgetvalue{/pgfplots/error bars/error mark}%
			{\pgfplotserrorbarsmark}%
			\pgfkeysgetvalue{/pgfplots/error bars/error mark options}%
			{\pgfplotserrorbarsmarkopts}%
			\path ##1 -- ##2;
		},
		/pgfplots/table/.cd,
		y index=5,
		y error expr={\thisrowno{3}-\thisrowno{5}},
		/pgfplots/box plot
	},
	box plot median/.style={
		/pgfplots/box plot
	}
}
\pgfplotsset{yticklabel style={text width=1.5em, align=right}}
\pgfplotsset{
  /pgfplots/xlabel near ticks/.style={
     /pgfplots/every axis x label/.style={
        at={(ticklabel cs:0.5)},anchor=near ticklabel
     }
  },
  /pgfplots/ylabel near ticks/.style={
     /pgfplots/every axis y label/.style={
        at={(ticklabel cs:0.5)},rotate=90,anchor=near ticklabel}
     }
 }

\pgfplotsset{compat=1.12}

\makeatletter

\usepackage{lineno,hyperref}

\journal{European Journal of Operational Research}

\biboptions{authoryear}

\begin{document}

\begin{frontmatter}

\title{Adaptive Solution Prediction for Combinatorial Optimization}


\author[addr1]{Yunzhuang Shen\corref{mycorrespondingauthor}}
\ead{s3640365@student.rmit.edu.au}
\cortext[mycorrespondingauthor]{Corresponding author}

\author[addr2]{Yuan Sun}
\ead{yuan.sun@latrobe.edu.au}

\author[addr1]{Xiaodong Li}
\ead{xiaodong.li@rmit.edu.au}

\author[addr3]{Andrew Eberhard}
\ead{andy.eberhard@rmit.edu.au}

\author[addr4]{Andreas Ernst}
\ead{andreas.ernst@monash.edu}

\address[addr1]{School of Computing Technologies, RMIT University, Melbourne, Australia}
\address[addr2]{La Trobe Business School, La Trobe University, Melbourne, Australia}
\address[addr3]{School of Science, RMIT University, Melbourne, Australia}
\address[addr4]{School of Mathematics, Monash University, Melbourne, Australia}

\begin{abstract}

This paper aims to predict optimal solutions for combinatorial optimization problems (COPs) via machine learning (ML). To find high-quality solutions efficiently, existing work uses a ML prediction of the optimal solution to guide heuristic search, where the ML model is trained offline under the supervision of solved problem instances with known optimal solutions. To predict the optimal solution with sufficient accuracy, it is critical to provide a ML model with adequate features that can effectively characterize decision variables. However, acquiring such features is challenging due to the high complexity of COPs. This paper proposes a framework that can better characterize decision variables by harnessing feedback from a heuristic search over several iterative steps, enabling an offline-trained ML model to predict the optimal solution in an adaptive manner. We refer to this approach as adaptive solution prediction (ASP). Specifically, we employ a set of statistical measures as features, which can extract useful information from feasible solutions found by a heuristic search and inform the ML model as to which value a decision variable is likely to take in high-quality solutions. Our experiments on three NP-hard COPs show that ASP substantially improves the prediction quality of an offline-trained ML model and achieves competitive results compared to several heuristic methods in terms of solution quality. Furthermore, we demonstrate that ASP can be used as a heuristic-pricing method for column generation, to boost an exact branch-and-price algorithm for solving the graph coloring problem.

\end{abstract}


\end{frontmatter}

\input{01-intro}
\input{02-approach}
\input{03-problems}
\input{04-result}

\input{05-cgbp}
\input{06-conclusion}
\bibliography{lib}

\end{document}


\begin{frontmatter}

\title{Adaptive Solution Prediction for Combinatorial Optimization Supplementary Materials}


\author[addr1]{Yunzhuang Shen\corref{mycorrespondingauthor}}
\ead{s3640365@student.rmit.edu.au}
\cortext[mycorrespondingauthor]{Corresponding author}

\author[addr2]{Yuan Sun}
\ead{yuan.sun@unimelb.edu.au}

\author[addr1]{Xiaodong Li}
\ead{xiaodong.li@rmit.edu.au}

\author[addr3]{Andrew Eberhard}
\ead{andy.eberhard@rmit.edu.au}

\author[addr4]{Andreas Ernst}
\ead{andreas.ernst@monash.edu}

\address[addr1]{School of Computing Technologies, RMIT University, Melbourne, Australia}
\address[addr2]{School of Computing and Information Systems, University of Melbourne, Melbourne, Australia}
\address[addr3]{School of Science, RMIT University, Melbourne, Australia}
\address[addr4]{School of Mathematics, Monash University, Melbourne, Australia}
\end{frontmatter}

\linenumbers

This manuscript reports additional experimental results. Specifically, Tables 1 and 2 contain results for comparing column generation using different heuristic methods on the set of $81$ graphs from the graph coloring benchmarks. Table 3 shows additional results for the branch-and-price method.

\begin{landscape}
\begin{table}[htb]
    \centering
    \caption{Numerical results of CG using different methods for solving the maximum weight independent set problem (MWISP). Each method is tested on $24$ problem instances generated from each graph under a cutoff time of $1800$ seconds. We report the total number of solved instances (i.e., RMPs), the average LP objective value of the RMPs, and the average solving time over the solved instances (if any).}
    \label{tab:cg-num-ret1}
    \resizebox{\linewidth}{!}{\begin{tabular}{@{}lrr|rrrrrrr|rrrrrrr|rrrrrrr@{}}
        \toprule
    \multirow{2}{*}{Graph} & \multirow{2}{*}{\# nodes} & \multirow{2}{*}{Density} & \multicolumn{7}{c}{\# optimally Solved} & \multicolumn{7}{c}{Solving time} & \multicolumn{7}{c}{LP objective}\\
    & & & ASP & SSSP & ACO & Gurobi & TSM & Fastwclq & LSCC & ASP & SSSP & ACO & Gurobi & TSM & Fastwclq & LSCC & ASP & SSSP & ACO & Gurobi & TSM & Fastwclq & LSCC \\
    \cmidrule(lr){1-3}\cmidrule(lr){4-10}\cmidrule(lr){11-17}\cmidrule(lr){18-24}

        wap08a & 1870 & 0.060 & 0 & 0 & 0 & 0 & 0 & 0 & 0 & N/A & N/A & N/A & N/A & N/A & N/A & N/A & 48.173 & \textbf{47.530} & 62.258 & 58.627 & 65.334 & 65.456 & 64.801 \\
        ash608GPIA & 1216 & 0.011 & \textbf{7} & 0 & 0 & 0 & 0 & 0 & 0 & \textbf{1560.3} & N/A & N/A & N/A & N/A & N/A & N/A & \textbf{3.333} & 3.372 & 3.368 & 3.394 & 3.394 & 3.385 & 3.394 \\
        wap07a & 1809 & 0.063 & 0 & 0 & 0 & 0 & 0 & 0 & 0 & N/A & N/A & N/A & N/A & N/A & N/A & N/A & 47.684 & \textbf{47.105} & 60.400 & 59.423 & 63.366 & 63.496 & 62.555 \\
        abb313GPIA & 1557 & 0.054 & 0 & 0 & 0 & 0 & 0 & 0 & 0 & N/A & N/A & N/A & N/A & N/A & N/A & N/A & \textbf{8.355} & 8.455 & 9.371 & 9.345 & 9.540 & 9.483 & 9.588 \\
        3-FullIns\_5 & 2030 & 0.016 & 0 & 0 & 0 & 0 & 0 & 0 & 0 & N/A & N/A & N/A & N/A & N/A & N/A & N/A & \textbf{5.609} & 5.697 & 8.763 & 7.684 & 9.672 & 6.607 & 9.349 \\
        DSJC1000.1 & 1000 & 0.099 & 0 & 0 & 0 & 0 & 0 & 0 & 0 & N/A & N/A & N/A & N/A & N/A & N/A & N/A & \textbf{19.871} & 19.984 & 21.028 & 22.107 & 22.296 & 20.837 & 21.645 \\
        3-Insertions\_5 & 1406 & 0.010 & 0 & 0 & 0 & 0 & 0 & 0 & 0 & N/A & N/A & N/A & N/A & N/A & N/A & N/A & \textbf{3.092} & 3.343 & 4.425 & 3.974 & 4.609 & 3.866 & 4.871 \\
        r1000.1 & 1000 & 0.029 & \textbf{24} & \textbf{24} & \textbf{24} & \textbf{24} & 11 & \textbf{24} & \textbf{24} & 182.7 & 210.3 & 263.0 & \textbf{145.4} & 987.7 & 633.1 & 300.2 & \textbf{20.000} & \textbf{20.000} & \textbf{20.000} & \textbf{20.000} & 20.017 & \textbf{20.000} & \textbf{20.000} \\
        DSJC1000.5 & 1000 & 0.500 & 0 & 0 & 0 & 0 & 0 & 0 & 0 & N/A & N/A & N/A & N/A & N/A & N/A & N/A & \textbf{77.825} & 78.028 & 78.401 & 90.751 & 88.568 & 84.175 & 87.938 \\
        flat1000\_60\_0 & 1000 & 0.492 & 0 & 0 & 0 & 0 & 0 & 0 & 0 & N/A & N/A & N/A & N/A & N/A & N/A & N/A & 76.014 & 76.308 & \textbf{75.755} & 89.235 & 87.034 & 81.522 & 85.709 \\
        flat1000\_76\_0 & 1000 & 0.494 & 0 & 0 & 0 & 0 & 0 & 0 & 0 & N/A & N/A & N/A & N/A & N/A & N/A & N/A & 77.178 & 77.293 & \textbf{77.134} & 89.580 & 87.434 & 82.189 & 86.803 \\
        flat1000\_50\_0 & 1000 & 0.490 & 0 & 0 & 0 & 0 & 0 & 0 & 0 & N/A & N/A & N/A & N/A & N/A & N/A & N/A & 53.562 & \textbf{51.760} & 61.696 & 88.794 & 86.420 & 66.964 & 80.768 \\
        5-FullIns\_4 & 1085 & 0.019 & \textbf{15} & 10 & 0 & 0 & 0 & 0 & 0 & \textbf{606.8} & 1228.4 & N/A & N/A & N/A & N/A & N/A & \textbf{7.283} & \textbf{7.283} & 7.616 & 7.407 & 7.715 & 7.325 & 7.902 \\
        will199GPIA & 701 & 0.029 & \textbf{24} & 23 & 0 & 5 & 0 & 0 & 0 & \textbf{171.7} & 774.3 & N/A & 1701.6 & N/A & N/A & N/A & \textbf{6.200} & \textbf{6.200} & 6.244 & 6.204 & 6.543 & 6.241 & 6.482 \\
        wap05a & 905 & 0.105 & \textbf{24} & \textbf{24} & 0 & \textbf{24} & \textbf{24} & 0 & 0 & 115.3 & \textbf{80.2} & N/A & 426.0 & 435.4 & N/A & N/A & \textbf{50.000} & \textbf{50.000} & 50.396 & \textbf{50.000} & \textbf{50.000} & 55.367 & 51.529 \\
        wap06a & 947 & 0.097 & \textbf{24} & 12 & 0 & 0 & 0 & 0 & 0 & \textbf{1080.8} & 1670.6 & N/A & N/A & N/A & N/A & N/A & \textbf{40.000} & 40.012 & 51.262 & 43.622 & 42.989 & 56.945 & 53.552 \\
        DSJC1000.9 & 1000 & 0.900 & \textbf{24} & \textbf{24} & \textbf{24} & 0 & \textbf{24} & 0 & 0 & \textbf{189.7} & 274.1 & 372.9 & N/A & 1579.4 & N/A & N/A & \textbf{214.855} & \textbf{214.855} & \textbf{214.855} & 228.839 & \textbf{214.855} & 214.973 & 219.215 \\
        DSJC500.1 & 500 & 0.100 & 0 & 0 & 0 & 0 & 0 & 0 & 0 & N/A & N/A & N/A & N/A & N/A & N/A & N/A & \textbf{11.341} & 11.415 & 12.141 & 12.888 & 13.394 & 11.363 & 11.848 \\
        2-FullIns\_5 & 852 & 0.034 & \textbf{23} & 0 & 0 & 0 & 0 & 0 & 0 & \textbf{610.2} & N/A & N/A & N/A & N/A & N/A & N/A & \textbf{4.708} & 4.737 & 5.610 & 5.247 & 5.593 & 4.973 & 6.105 \\
        4-Insertions\_4 & 475 & 0.016 & 0 & 0 & 0 & 0 & 0 & 0 & 0 & N/A & N/A & N/A & N/A & N/A & N/A & N/A & \textbf{2.465} & 2.592 & 2.745 & 2.852 & 3.124 & 2.641 & 3.025 \\
        2-Insertions\_5 & 597 & 0.022 & 0 & 0 & 0 & 0 & 0 & 0 & 0 & N/A & N/A & N/A & N/A & N/A & N/A & N/A & \textbf{2.820} & 2.920 & 3.458 & 3.349 & 3.730 & 2.978 & 3.747 \\
        4-FullIns\_4 & 690 & 0.028 & \textbf{24} & \textbf{24} & 0 & \textbf{24} & 0 & 22 & 0 & \textbf{395.7} & 636.0 & N/A & 778.0 & N/A & 1034.5 & N/A & \textbf{6.329} & \textbf{6.329} & 6.389 & \textbf{6.329} & 6.432 & \textbf{6.329} & 6.575 \\
        r1000.5 & 1000 & 0.477 & \textbf{24} & 14 & \textbf{24} & 0 & \textbf{24} & 0 & 0 & \textbf{411.3} & 1421.3 & 1155.1 & N/A & 786.8 & N/A & N/A & \textbf{234.000} & 234.053 & \textbf{234.000} & 268.258 & \textbf{234.000} & 249.398 & 243.494 \\
        DSJC500.5 & 500 & 0.502 & 0 & 0 & 0 & 0 & 0 & 0 & 0 & N/A & N/A & N/A & N/A & N/A & N/A & N/A & 42.472 & 42.453 & 42.628 & 50.549 & \textbf{42.347} & 43.559 & 44.098 \\
        1-Insertions\_6 & 607 & 0.034 & 0 & 0 & 0 & 0 & 0 & 0 & 0 & N/A & N/A & N/A & N/A & N/A & N/A & N/A & \textbf{3.188} & 3.227 & 4.264 & 3.918 & 4.029 & 3.324 & 4.664 \\
        le450\_5a & 450 & 0.057 & 0 & 0 & 0 & 0 & 0 & 0 & 0 & N/A & N/A & N/A & N/A & N/A & N/A & N/A & \textbf{5.991} & 6.247 & 6.031 & 7.017 & 8.880 & 6.129 & 7.336 \\
        le450\_5b & 450 & 0.057 & 0 & 0 & 0 & 0 & 0 & 0 & 0 & N/A & N/A & N/A & N/A & N/A & N/A & N/A & 6.050 & 6.299 & \textbf{5.900} & 7.025 & 8.837 & 6.113 & 7.321 \\
        r1000.1c & 1000 & 0.971 & \textbf{24} & \textbf{24} & \textbf{24} & 0 & \textbf{24} & \textbf{24} & 1 & \textbf{19.2} & 21.5 & 128.8 & N/A & 68.7 & 53.2 & 1774.6 & \textbf{95.057} & \textbf{95.057} & \textbf{95.057} & 95.197 & \textbf{95.057} & \textbf{95.057} & 95.071 \\
        le450\_25a & 450 & 0.082 & 24 & 24 & 24 & 24 & 24 & 24 & 24 & 10.7 & 11.4 & 90.8 & 11.5 & \textbf{10.5} & 637.3 & 98.3 & 25.000 & 25.000 & 25.000 & 25.000 & 25.000 & 25.000 & 25.000 \\
        le450\_15d & 450 & 0.166 & 0 & 0 & 0 & 0 & 0 & 0 & 0 & N/A & N/A & N/A & N/A & N/A & N/A & N/A & \textbf{17.262} & 17.725 & 21.902 & 20.015 & 21.687 & 21.862 & 21.212 \\
        le450\_15c & 450 & 0.165 & 0 & 0 & 0 & 0 & 0 & 0 & 0 & N/A & N/A & N/A & N/A & N/A & N/A & N/A & \textbf{17.122} & 17.553 & 21.913 & 20.061 & 21.458 & 21.734 & 21.002 \\
        le450\_15b & 450 & 0.081 & \textbf{24} & \textbf{24} & 0 & \textbf{24} & 0 & 0 & 0 & \textbf{85.3} & 109.4 & N/A & 1265.6 & N/A & N/A & N/A & \textbf{15.000} & \textbf{15.000} & 17.082 & \textbf{15.000} & 15.932 & 17.327 & 15.921 \\
        le450\_15a & 450 & 0.081 & \textbf{24} & \textbf{24} & 0 & \textbf{24} & 0 & 0 & 0 & \textbf{89.6} & 127.6 & N/A & 1319.0 & N/A & N/A & N/A & \textbf{15.000} & \textbf{15.000} & 17.186 & \textbf{15.000} & 15.870 & 17.557 & 15.941 \\
        le450\_25c & 450 & 0.172 & \textbf{24} & \textbf{24} & 0 & 0 & 0 & 0 & 0 & \textbf{178.5} & 254.0 & N/A & N/A & N/A & N/A & N/A & \textbf{25.000} & \textbf{25.000} & 28.468 & 25.428 & 25.942 & 29.577 & 27.268 \\
        le450\_25d & 450 & 0.172 & \textbf{24} & \textbf{24} & 0 & 0 & 0 & 0 & 0 & \textbf{171.1} & 188.5 & N/A & N/A & N/A & N/A & N/A & \textbf{25.000} & \textbf{25.000} & 28.011 & 25.408 & 26.091 & 29.055 & 27.246 \\
        le450\_5c & 450 & 0.097 & 0 & 0 & \textbf{7} & 0 & 0 & 0 & 0 & N/A & N/A & \textbf{589.3} & N/A & N/A & N/A & N/A & 5.231 & 5.346 & \textbf{5.109} & 7.423 & 7.477 & 5.440 & 7.071 \\
        queen16\_16 & 256 & 0.387 & \textbf{24} & \textbf{24} & \textbf{24} & \textbf{24} & 0 & 17 & \textbf{24} & \textbf{29.1} & 34.8 & 135.0 & 418.2 & N/A & 419.5 & 697.8 & \textbf{16.000} & \textbf{16.000} & \textbf{16.000} & \textbf{16.000} & 16.805 & 16.070 & \textbf{16.000} \\
        le450\_5d & 450 & 0.097 & 1 & 0 & \textbf{11} & 0 & 0 & 0 & 0 & 1637.2 & N/A & \textbf{713.6} & N/A & N/A & N/A & N/A & 5.143 & 5.343 & \textbf{5.092} & 7.395 & 7.565 & 5.351 & 7.037 \\
        DSJC500.9 & 500 & 0.901 & \textbf{24} & \textbf{24} & \textbf{24} & 0 & \textbf{24} & \textbf{24} & \textbf{24} & \textbf{10.4} & 12.6 & 168.4 & N/A & 71.7 & 111.6 & 389.0 & \textbf{122.306} & \textbf{122.306} & \textbf{122.306} & 125.304 & \textbf{122.306} & \textbf{122.306} & \textbf{122.306} \\
        3-FullIns\_4 & 405 & 0.043 & \textbf{24} & \textbf{24} & \textbf{24} & \textbf{24} & \textbf{24} & \textbf{24} & 23 & \textbf{30.8} & 49.4 & 533.3 & 119.7 & 327.7 & 220.2 & 656.9 & 5.392 & 5.392 & 5.392 & 5.392 & 5.392 & 5.392 & 5.392 \\
        \bottomrule
    \end{tabular}}
\end{table}
\end{landscape}

\begin{landscape}
\begin{table}[htb]
    \centering
    \caption{Numerical results of CG using different methods for solving MWISP (Table~\ref{tab:cg-num-ret1} continued).}
    \label{tab:cg-num-ret2}
    
    \resizebox{\linewidth}{!}{\begin{tabular}{@{}lrr|rrrrrrr|rrrrrrr|rrrrrrr@{}}
        \toprule
            \multirow{2}{*}{Graph} & \multirow{2}{*}{\# nodes} & \multirow{2}{*}{Density} & \multicolumn{7}{c}{\# optimally Solved} & \multicolumn{7}{c}{Solving time} & \multicolumn{7}{c}{LP objective}\\
    & & & ASP & SSSP & ACO & Gurobi & TSM & Fastwclq & LSCC & ASP & SSSP & ACO & Gurobi & TSM & Fastwclq & LSCC & ASP & SSSP & ACO & Gurobi & TSM & Fastwclq & LSCC \\
    \cmidrule(lr){1-3}\cmidrule(lr){4-10}\cmidrule(lr){11-17}\cmidrule(lr){18-24}
        queen15\_15 & 225 & 0.411 & \textbf{24} & \textbf{24} & \textbf{24} & \textbf{24} & 0 & 22 & \textbf{24} & \textbf{16.9} & 18.9 & 58.1 & 269.0 & N/A & 184.3 & 373.8 & \textbf{15.000} & \textbf{15.000} & \textbf{15.000} & \textbf{15.000} & 15.326 & 15.092 & \textbf{15.000} \\
        3-Insertions\_4 & 281 & 0.027 & 0 & 0 & 0 & 0 & 0 & 0 & 0 & N/A & N/A & N/A & N/A & N/A & N/A & N/A & \textbf{2.463} & 2.533 & 2.619 & 2.604 & 2.790 & 2.477 & 2.796 \\
        school1 & 385 & 0.258 & 0 & 0 & 0 & 0 & 0 & 0 & 0 & N/A & N/A & N/A & N/A & N/A & N/A & N/A & 15.914 & \textbf{15.420} & 21.960 & 19.523 & 15.870 & 16.188 & 18.333 \\
        DSJR500.5 & 500 & 0.472 & \textbf{24} & \textbf{24} & \textbf{24} & 0 & \textbf{24} & 21 & \textbf{24} & \textbf{10.3} & 16.4 & 159.8 & N/A & 38.9 & 517.4 & 274.4 & \textbf{122.000} & \textbf{122.000} & \textbf{122.000} & 132.687 & \textbf{122.000} & 122.399 & \textbf{122.000} \\
        flat300\_26\_0 & 300 & 0.482 & 0 & 0 & 0 & 0 & 0 & 0 & 0 & N/A & N/A & N/A & N/A & N/A & N/A & N/A & \textbf{26.000} & \textbf{26.000} & \textbf{26.000} & 30.654 & \textbf{26.000} & 26.254 & \textbf{26.000} \\
        flat300\_28\_0 & 300 & 0.484 & 0 & 0 & 9 & 0 & \textbf{24} & 19 & 3 & N/A & N/A & 1465.2 & N/A & \textbf{225.0} & 1046.9 & 1487.7 & 27.527 & 27.525 & \textbf{27.520} & 30.780 & \textbf{27.520} & 27.698 & \textbf{27.520} \\
        school1\_nsh & 352 & 0.237 & 0 & 0 & 0 & 0 & 0 & 0 & 0 & N/A & N/A & N/A & N/A & N/A & N/A & N/A & 15.396 & 15.657 & 21.559 & 18.142 & \textbf{15.176} & 17.250 & 17.854 \\
        flat300\_20\_0 & 300 & 0.477 & 0 & 0 & 0 & 0 & 0 & 0 & 0 & N/A & N/A & N/A & N/A & N/A & N/A & N/A & \textbf{20.000} & \textbf{20.000} & \textbf{20.000} & 24.526 & 20.522 & \textbf{20.000} & 20.896 \\
        DSJC250.1 & 250 & 0.103 & 0 & 0 & 0 & 0 & 0 & 0 & 0 & N/A & N/A & N/A & N/A & N/A & N/A & N/A & 7.006 & 7.079 & 7.625 & 7.466 & 8.033 & \textbf{7.004} & 7.330 \\
        queen14\_14 & 196 & 0.438 & \textbf{24} & \textbf{24} & \textbf{24} & \textbf{24} & \textbf{24} & 23 & \textbf{24} & \textbf{10.5} & 12.1 & 31.2 & 183.7 & 220.1 & 166.9 & 218.8 & \textbf{14.000} & \textbf{14.000} & \textbf{14.000} & \textbf{14.000} & \textbf{14.000} & 14.017 & \textbf{14.000} \\
        DSJC250.5 & 250 & 0.503 & 14 & 5 & \textbf{24} & 0 & \textbf{24} & \textbf{24} & 17 & 1302.8 & 1539.4 & 310.4 & N/A & \textbf{83.9} & 302.2 & 1009.6 & \textbf{25.165} & 25.166 & \textbf{25.165} & 26.714 & \textbf{25.165} & \textbf{25.165} & \textbf{25.165} \\
        queen13\_13 & 169 & 0.469 & 24 & 24 & 24 & 24 & 24 & 24 & 24 & \textbf{4.9} & 5.5 & 13.4 & 91.3 & 42.9 & 47.0 & 104.8 & 13.000 & 13.000 & 13.000 & 13.000 & 13.000 & 13.000 & 13.000 \\
        DSJR500.1c & 500 & 0.972 & \textbf{24} & \textbf{24} & \textbf{24} & 11 & \textbf{24} & \textbf{24} & \textbf{24} & 1.5 & 1.8 & 62.2 & 381.0 & \textbf{1.2} & 1.4 & 26.5 & \textbf{84.136} & \textbf{84.136} & \textbf{84.136} & 84.158 & \textbf{84.136} & \textbf{84.136} & \textbf{84.136} \\
        1-FullIns\_5 & 282 & 0.082 & \textbf{24} & \textbf{24} & 1 & 15 & 0 & \textbf{24} & 0 & 13.7 & \textbf{13.0} & 1041.8 & 1448.3 & N/A & 118.7 & N/A & \textbf{3.909} & \textbf{3.909} & 3.957 & 3.910 & 3.939 & \textbf{3.909} & 3.974 \\
        1-Insertions\_5 & 202 & 0.060 & \textbf{24} & 0 & 0 & 0 & 0 & \textbf{24} & 0 & 667.2 & N/A & N/A & N/A & N/A & \textbf{449.8} & N/A & \textbf{2.943} & 2.953 & 3.130 & 2.963 & 2.988 & \textbf{2.943} & 3.149 \\
        queen12\_12 & 144 & 0.504 & 24 & 24 & 24 & 24 & 24 & 24 & 24 & \textbf{3.2} & 3.4 & 7.2 & 60.2 & 13.3 & 29.1 & 60.3 & 12.000 & 12.000 & 12.000 & 12.000 & 12.000 & 12.000 & 12.000 \\
        2-FullIns\_4 & 212 & 0.072 & 24 & 24 & 24 & 24 & 24 & 24 & 24 & \textbf{4.0} & \textbf{4.0} & 157.3 & 14.0 & 15.6 & 17.9 & 298.0 & 4.485 & 4.485 & 4.485 & 4.485 & 4.485 & 4.485 & 4.485 \\
        DSJC250.9 & 250 & 0.896 & \textbf{24} & \textbf{24} & \textbf{24} & 0 & \textbf{24} & \textbf{24} & \textbf{24} & \textbf{1.4} & 1.5 & 6.2 & N/A & 4.9 & 6.5 & 16.4 & \textbf{70.392} & \textbf{70.392} & \textbf{70.392} & 71.240 & \textbf{70.392} & \textbf{70.392} & \textbf{70.392} \\
        2-Insertions\_4 & 149 & 0.049 & \textbf{24} & 0 & 0 & \textbf{24} & \textbf{24} & \textbf{24} & 0 & 466.7 & N/A & N/A & 439.7 & 1328.3 & \textbf{239.0} & N/A & \textbf{2.560} & 2.567 & 2.592 & \textbf{2.560} & \textbf{2.560} & \textbf{2.560} & 2.600 \\
        r250.5 & 250 & 0.477 & 24 & 24 & 24 & 24 & 24 & 24 & 24 & \textbf{1.5} & 2.2 & 11.3 & 476.3 & 4.3 & 61.4 & 26.1 & 65.000 & 65.000 & 65.000 & 65.000 & 65.000 & 65.000 & 65.000 \\
        queen11\_11 & 121 & 0.545 & 24 & 24 & 24 & 24 & 24 & 24 & 24 & \textbf{1.8} & 2.0 & 3.6 & 36.3 & 5.9 & 14.5 & 32.0 & 11.000 & 11.000 & 11.000 & 11.000 & 11.000 & 11.000 & 11.000 \\
        mug100\_1 & 100 & 0.034 & 24 & 24 & 24 & 24 & 24 & 24 & 24 & \textbf{4.0} & 9.5 & 27.5 & 17.5 & 18.6 & 54.8 & 48.3 & 3.030 & 3.030 & 3.030 & 3.030 & 3.030 & 3.030 & 3.030 \\
        mug100\_25 & 100 & 0.034 & 24 & 24 & 24 & 24 & 24 & 24 & 24 & \textbf{3.2} & 6.6 & 21.6 & 16.4 & 15.8 & 48.8 & 42.6 & 3.030 & 3.030 & 3.030 & 3.030 & 3.030 & 3.030 & 3.030 \\
        mug88\_1 & 88 & 0.038 & 24 & 24 & 24 & 24 & 24 & 24 & 24 & \textbf{1.7} & 3.4 & 9.3 & 8.7 & 6.5 & 25.7 & 18.8 & 3.034 & 3.034 & 3.034 & 3.034 & 3.034 & 3.034 & 3.034 \\
        myciel7 & 191 & 0.130 & \textbf{24} & \textbf{24} & 0 & \textbf{24} & 9 & \textbf{24} & 0 & 127.4 & 132.8 & N/A & 1097.1 & 1330.7 & \textbf{124.2} & N/A & \textbf{4.095} & \textbf{4.095} & 4.160 & \textbf{4.095} & 4.113 & \textbf{4.095} & 4.097 \\
        DSJC125.1 & 125 & 0.095 & 0 & 0 & 0 & \textbf{24} & \textbf{24} & 0 & 0 & N/A & N/A & N/A & 495.2 & \textbf{75.2} & N/A & N/A & \textbf{4.454} & 4.458 & 4.565 & \textbf{4.454} & \textbf{4.454} & 4.484 & 4.461 \\
        mug88\_25 & 88 & 0.038 & 24 & 24 & 24 & 24 & 24 & 24 & 24 & \textbf{1.9} & 3.9 & 9.3 & 9.4 & 8.1 & 23.1 & 22.6 & 3.034 & 3.034 & 3.034 & 3.034 & 3.034 & 3.034 & 3.034 \\
        queen10\_10 & 100 & 0.594 & 24 & 24 & 24 & 24 & 24 & 24 & 24 & \textbf{1.1} & 1.2 & 1.7 & 28.1 & 3.7 & 8.4 & 21.1 & 10.000 & 10.000 & 10.000 & 10.000 & 10.000 & 10.000 & 10.000 \\
        DSJC125.5 & 125 & 0.502 & 24 & 24 & 24 & 24 & 24 & 24 & 24 & 12.0 & 56.2 & 7.8 & 908.7 & \textbf{4.8} & 14.3 & 53.5 & 15.727 & 15.727 & 15.727 & 15.727 & 15.727 & 15.727 & 15.727 \\
4-Insertions\_3 & 79 & 0.051 & 24 & 24 & 24 & 24 & 24 & 24 & 24 & 9.6 & 362.8 & 299.8 & 16.7 & \textbf{8.8} & 25.6 & 61.5 & 2.276 & 2.276 & 2.276 & 2.276 & 2.276 & 2.276 & 2.276 \\
        queen9\_9 & 81 & 0.652 & 24 & 24 & 24 & 24 & 24 & 24 & 24 & \textbf{0.5} & \textbf{0.5} & 0.7 & 9.8 & 1.4 & 3.5 & 7.0 & 9.000 & 9.000 & 9.000 & 9.000 & 9.000 & 9.000 & 9.000 \\
        r125.5 & 125 & 0.495 & 24 & 24 & 24 & 24 & 24 & 24 & 24 & \textbf{0.2} & 0.3 & 0.7 & 3.5 & 0.3 & 1.0 & 1.2 & 36.000 & 36.000 & 36.000 & 36.000 & 36.000 & 36.000 & 36.000 \\
        DSJC125.9 & 125 & 0.898 & 24 & 24 & 24 & 24 & 24 & 24 & 24 & \textbf{0.2} & \textbf{0.2} & 0.5 & 53.7 & 0.3 & 0.4 & 0.5 & 42.727 & 42.727 & 42.727 & 42.727 & 42.727 & 42.727 & 42.727 \\
        1-FullIns\_4 & 93 & 0.139 & 24 & 24 & 24 & 24 & 24 & 24 & 24 & 0.3 & \textbf{0.2} & 1.3 & 2.0 & 0.9 & 0.7 & 1.4 & 3.633 & 3.633 & 3.633 & 3.633 & 3.633 & 3.633 & 3.633 \\
        1-Insertions\_4 & 67 & 0.105 & 24 & 24 & 24 & 24 & 24 & 24 & 24 & \textbf{1.4} & 3.3 & 39.0 & 9.2 & 2.4 & 4.5 & 19.2 & 2.774 & 2.774 & 2.774 & 2.774 & 2.774 & 2.774 & 2.774 \\
        myciel6 & 95 & 0.169 & 24 & 24 & 24 & 24 & 24 & 24 & 24 & 2.1 & \textbf{1.9} & 518.7 & 28.6 & 12.1 & 9.9 & 76.1 & 3.834 & 3.834 & 3.834 & 3.834 & 3.834 & 3.834 & 3.834 \\
        3-Insertions\_3 & 56 & 0.071 & 24 & 24 & 24 & 24 & 24 & 24 & 24 & 1.2 & 6.8 & 13.2 & 3.5 & \textbf{1.1} & 4.1 & 11.2 & 2.334 & 2.334 & 2.334 & 2.334 & 2.334 & 2.334 & 2.334 \\
        queen8\_8 & 64 & 0.722 & 24 & 24 & 24 & 24 & 24 & 24 & 24 & \textbf{0.1} & \textbf{0.1} & \textbf{0.1} & 1.3 & 0.2 & 0.4 & 0.4 & 8.444 & 8.444 & 8.444 & 8.444 & 8.444 & 8.444 & 8.444 \\
        2-Insertions\_3 & 37 & 0.108 & 24 & 24 & 24 & 24 & 24 & 24 & 24 & \textbf{0.1} & 0.2 & 0.3 & 0.6 & 0.2 & 0.4 & 0.4 & 2.423 & 2.423 & 2.423 & 2.423 & 2.423 & 2.423 & 2.423 \\
        myciel5 & 47 & 0.218 & 24 & 24 & 24 & 24 & 24 & 24 & 24 & \textbf{0.1} & \textbf{0.1} & 1.2 & 1.1 & 0.3 & 0.7 & 1.6 & 3.553 & 3.553 & 3.553 & 3.553 & 3.553 & 3.553 & 3.553 \\
        myciel4 & 23 & 0.281 & 24 & 24 & 24 & 24 & 24 & 24 & 24 & \textbf{0.0} & \textbf{0.0} & \textbf{0.0} & 0.1 & \textbf{0.0} & \textbf{0.0} & \textbf{0.0} & 3.245 & 3.245 & 3.245 & 3.245 & 3.245 & 3.245 & 3.245 \\
    \cmidrule(lr){1-3}\cmidrule(lr){4-10}\cmidrule(lr){11-17}\cmidrule(lr){18-24}
        statistics &  &  & \textbf{1212} & 1096 & 916 & 895 & 932 & 964 & 812 & \textbf{12.879} & 12.934 & 13.615 & 13.838 & 13.903 & 13.494 & 13.927 & N/A & N/A & N/A & N/A & N/A & N/A\\
        \bottomrule

    \end{tabular}}
\end{table}
\end{landscape}

\begin{table}[htb]
    \centering
    \caption{The rest of the results for graphs where B\&P with different pricing methods have comparable performances. Each method performs $24$ seeded runs for a graph. The solving time is an average over optimally solved GCPs, and `N/A' is placed if a method does not solve any instances on that graph. Similarly, the duality gap for a graph is an average over instances whose root LP has been solved to optimality, and `N/A' is placed if the root LP cannot be solved in any run.}
    \label{tab:bp2}
    \resizebox{\linewidth}{!}{\begin{tabular}{@{}lrr|rrr|rrr|rrr|rrr@{}}
        \toprule
        \multirow{2}{*}{Graph} & \multirow{2}{*}{\# nodes} & \multirow{2}{*}{Density} & \multicolumn{3}{c}{\# optimally Solved} & \multicolumn{3}{c}{Solving time}  & \multicolumn{3}{c}{Duality gap} & \multicolumn{3}{c}{\# root LP solved} \\
    & & & ASP & SSSP & Greedy & ASP & SSSP & Greedy & ASP & SSSP & Greedy & ASP & SSSP & Greedy \\
    \cmidrule(lr){1-3}\cmidrule(lr){4-6}\cmidrule(lr){7-9}\cmidrule(lr){10-12}\cmidrule(lr){13-15}
     r250.5 & 235 & 0.508 & 24 & 24 & 24 & 44.9 & 58.8 & 46.6 & 0.000 & 0.000 & 0.000 & 24 & 24 & 24 \\
    DSJR500.1c & 311 & 0.972 & 24 & 24 & 24 & 99.9 & 99.7 & 109.6 & 0.000 & 0.000 & 0.000 & 24 & 24 & 24 \\
    r125.5 & 109 & 0.565 & 24 & 24 & 24 & 12.7 & 12.1 & 13.3 & 0.000 & 0.000 & 0.000 & 24 & 24 & 24 \\
  DSJC125.9 & 125 & 0.898 & 24 & 24 & 24 & 27.5 & 32.0 & 29.1 & 0.000 & 0.000 & 0.000 & 24 & 24 & 24 \\
   le450\_25a & 264 & 0.336 & 24 & 24 & 24 & 13.8 & 13.6 & 13.1 & 0.000 & 0.000 & 0.000 & 24 & 24 & 24 \\
    qg.order30 & 900 & 0.129 & 24 & 24 & 24 & 60.7 & 58.6 & 60.6 & 0.000 & 0.000 & 0.000 & 24 & 24 & 24 \\
        qg.order60 & 3600 & 0.066 & 24 & 24 & 24 & 678.0 & 631.7 & 634.1 & 0.000 & 0.000 & 0.000 & 24 & 24 & 24 \\
        wap05a & 665 & 0.314 & 24 & 24 & 24 & 76.3 & 74.6 & 76.5 & 0.000 & 0.000 & 0.000 & 24 & 24 & 24 \\
        le450\_25b & 294 & 0.290 & 24 & 24 & 24 & 14.8 & 14.7 & 16.5 & 0.000 & 0.000 & 0.000 & 24 & 24 & 24 \\
  qg.order100 & 10000 & 0.040 & 24 & 24 & 24 & 6528.6 & 6259.1 & 6127.9 & 0.000 & 0.000 & 0.000 & 24 & 24 & 24 \\
   queen10\_10 & 100 & 0.594 & 0 & 0 & 0 & N/A & N/A & N/A & 16.667 & 16.351 & 16.035 & 24 & 24 & 24 \\
    DSJC125.1 & 125 & 0.095 & 0 & 0 & 0 & N/A & N/A & N/A & 16.667 & 16.667 & 16.667 & 24 & 24 & 24 \\

    3-Insertions\_4 & 281 & 0.053 & 0 & 0 & 0 & N/A & N/A & N/A & 40.000 & 40.000 & 40.000 & 24 & 24 & 24 \\
    5-FullIns\_4 & 43 & 0.651 & 0 & 0 & 0 & N/A & N/A & N/A & 11.111 & 11.111 & 11.111 & 24 & 24 & 24 \\
    queen13\_13 & 169 & 0.469 & 0 & 0 & 0 & N/A & N/A & N/A & 13.333 & 13.333 & 13.333 & 24 & 24 & 24 \\
    4-FullIns\_5 & 113 & 0.293 & 0 & 0 & 0 & N/A & N/A & N/A & 22.222 & 22.222 & 22.222 & 24 & 24 & 24 \\
    3-Insertions\_3 & 56 & 0.143 & 0 & 0 & 0 & N/A & N/A & N/A & 25.000 & 25.000 & 25.000 & 24 & 24 & 24 \\
    DSJC500.9 & 500 & 0.901 & 0 & 0 & 0 & N/A & N/A & N/A & 6.818 & 6.818 & 6.818 & 24 & 24 & 24 \\
    DSJC250.5 & 250 & 0.503 & 0 & 0 & 0 & N/A & N/A & N/A & 16.129 & 16.129 & 16.129 & 24 & 24 & 24 \\
    2-Insertions\_4 & 149 & 0.098 & 0 & 0 & 0 & N/A & N/A & N/A & 40.000 & 40.000 & 40.000 & 24 & 24 & 24 \\
    queen14\_14 & 196 & 0.438 & 0 & 0 & 0 & N/A & N/A & N/A & 12.500 & 12.500 & 12.500 & 24 & 24 & 24 \\
    queen11\_11 & 121 & 0.545 & 0 & 0 & 0 & N/A & N/A & N/A & 15.385 & 15.385 & 15.385 & 24 & 24 & 24 \\
    latin\_square\_10 & 900 & 0.760 & 0 & 0 & 0 & N/A & N/A & N/A & 17.431 & 17.431 & 17.431 & 24 & 24 & 24 \\
    myciel7 & 191 & 0.260 & 0 & 0 & 0 & N/A & N/A & N/A & 37.500 & 37.500 & 37.500 & 24 & 24 & 24 \\
    queen12\_12 & 144 & 0.504 & 0 & 0 & 0 & N/A & N/A & N/A & 14.286 & 14.286 & 14.286 & 24 & 24 & 24 \\
    2-FullIns\_4 & 41 & 0.393 & 0 & 0 & 0 & N/A & N/A & N/A & 16.667 & 16.667 & 16.667 & 24 & 24 & 24 \\
    4-FullIns\_4 & 37 & 0.649 & 0 & 0 & 0 & N/A & N/A & N/A & 12.500 & 12.500 & 12.500 & 24 & 24 & 24 \\
    flat300\_28\_0 & 300 & 0.967 & 0 & 0 & 0 & N/A & N/A & N/A & 17.647 & 17.647 & 17.647 & 24 & 24 & 24 \\
    2-FullIns\_5 & 92 & 0.263 & 0 & 0 & 0 & N/A & N/A & N/A & 28.571 & 28.571 & 28.571 & 24 & 24 & 24 \\
    3-FullIns\_4 & 43 & 0.441 & 0 & 0 & 0 & N/A & N/A & N/A & 14.286 & 14.286 & 14.286 & 24 & 24 & 24 \\
    3-FullIns\_5 & 103 & 0.276 & 0 & 0 & 0 & N/A & N/A & N/A & 25.000 & 25.000 & 25.000 & 24 & 24 & 24 \\
    4-Insertions\_3 & 79 & 0.101 & 0 & 0 & 0 & N/A & N/A & N/A & 25.000 & 25.000 & 25.000 & 24 & 24 & 24 \\
        \bottomrule
    \end{tabular}}
\end{table}

%% file: 01-intro.tex
\section{Introduction}

Many real-world problems can be formulated as combinatorial optimization problems (COPs), which can be expressed by minimizing a linear objective function of integer decision variables, subject to a set of linear inequalities~\citep{wolsey1999integer}. COPs are generally NP-hard to solve, and many effects have been made in past decades in devising heuristic methods to find good solutions efficiently. Among these heuristics, the most generic type is the one based on mixed-integer-programming (MIP) techniques~\citep{fischetti2010heuristics} and can be readily applied to general COPs. On the other hand, for some well-studied problems, specialized methods have been developed by human experts by exploiting the structure of a specific problem, such as Lin-Kernighan-Helsgaun~\citep{helsgaun2017extension} for the traveling salesman problem (TSP). Somewhere in between the two extremes are the metaheuristics~\citep{glover2006handbook}, which are based on certain general assumptions and may be adapted to specific problems to some extent. An example is the class of estimation-of-distribution algorithms (EDAs)~\citep{hauschild2011introduction}, which alternates between sampling solutions from the search space and learning a probabilistic model over good samples to improve the sample efficiency, here referred to as \textit{online} learning. Meanwhile, exact optimization methods have also been improved substantially, including generic MIP solvers~\citep{50year2010, miplib2017} such as Gurobi~\citep{gurobi2018gurobi} and specialized ones such as Concorde~\citep{applegate2011traveling} for TSP.

In recent years, leveraging machine learning (ML) to develop heuristic methods has attracted a lot of attention~\citep{bengio2020machine,rl2021nina}. Compared with conventional methods, ML can automatically learn knowledge \textit{offline} from historical data of a COP and apply it in solving unseen but similar problem instances. In particular, heuristic methods based on solution prediction~\citep{li2018combinatorial,lauri2019fine,sun2019using,ding2020accelerating,zhang2020nlocalsat,abbasi2020predicting} have been developed: a ML model is trained offline under the supervision of solved problem instances with known optimal solutions; given an unseen instance, the offline-trained ML model is used to make a prediction of its optimal solution to guide a heuristic search.

To predict the optimal solution with sufficient accuracy, it is critical to provide a ML model with adequate features that can characterize decision variables effectively. However, acquiring such features is challenging due to the high complexity of COPs, particularly for problems that are highly constrained and/or contain many local-optimal solutions. Moreover, with the growth of the problem size, the same set of features can be less sufficient to represent and distinguish the decision variables. As a result, the performance of a solution-prediction-based method can decrease significantly for solving large-scale problems~\citep{joshi2019learning}. Existing methods typically guide the heuristic search using the ML prediction made only once \textit{prior to the search process}. This type of method can be inefficient when the ML model fails to make a sufficiently accurate prediction due to inadequate features.

\begin{figure}[ht!]
	\centering
	\includegraphics[width=\textwidth]{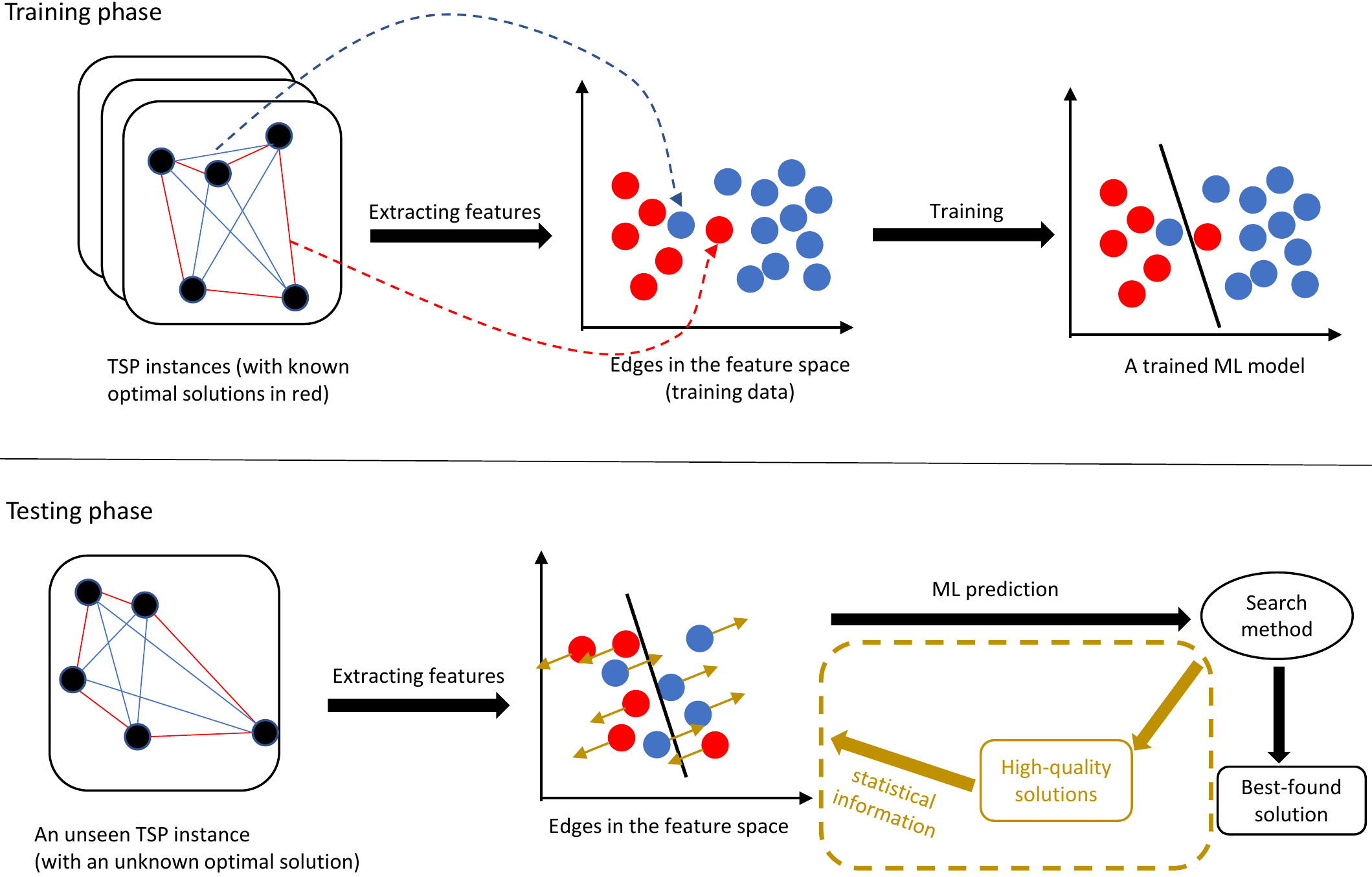}
	\caption{The adaptive solution prediction framework, illustrated on the TSP. In the training phase, some TSP instances are solved to optimality in advance; using these instances, we can construct a training data set where each training example (denoted by a dot) corresponds to an edge in a TSP instance; a ML model is then trained to best separate the optimal and non-optimal edges. In the testing phase, given an unseen TSP instance with its edges represented in the feature space, we can use the offline-trained ML model to predict the optimal solution values of these edges and use the ML prediction to guide a search method for finding high-quality solutions. We leverage feedback from the search to refine the features of edges such that the optimal and non-optimal edges are more separable, as indicated by the golden arrows. Specifically, from an evolving set of best-found solutions, we extract more accurate statistical information and use it as a part of the feature representation for edges. It is important to note that the ML model is not retrained in the testing phase. The diagram without the feedback mechanism describes the existing solution-prediction methods.}
	\label{fig:asp}
\end{figure}

This paper aims to improve the accuracy of the ML prediction to boost heuristic search. We propose a framework that can better characterize decision variables by leveraging useful information progressively collected \textit{during the search process}, enabling an offline-trained ML model to predict the optimal solution in an adaptive manner. We refer to our method as adaptive solution prediction (ASP), as illustrated in Figure~\ref{fig:asp}. Our key innovation is the feedback mechanism highlighted in the golden dashed box. Specifically, we employ a set of statistical measures as part of the feature representation of a decision variable (e.g., an edge in the TSP); these statistical features can extract useful information from a pool of feasible solutions and suggest which value the variable is likely to take in high-quality solutions. During the search process, we recompute the values of these statistical features for every decision variable. More accurate statistical information can be obtained to better characterize these variables, as the heuristic search continuously supplies better-quality solutions into the solution pool. As a result, the offline-trained ML model can better predict the optimal solution, which in turn better guides the search. The feedback mechanism of ASP can be viewed as a form of online learning similar to EDAs, which will be further discussed in Section~\ref{subsec:MHcompare} after presenting our method. Since ASP benefits from both offline and online learning, it can generalize better than existing ML-based methods and is more efficient than conventional EDAs.

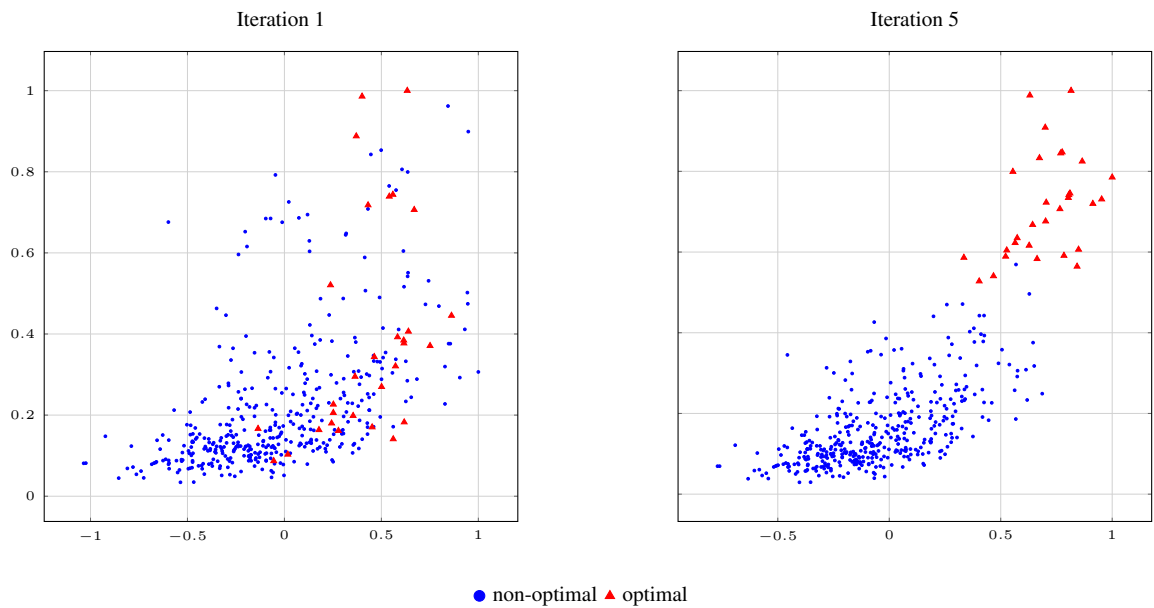
\begin{figure}[ht!]
    \centering
    \begin{tikzpicture}[baseline={(0,0)}]
    \begin{groupplot}[group style = {group size = 2 by 1, horizontal sep = 60pt},height=0.475\textwidth, width=0.475\textwidth, grid style={line width=.1pt, draw=gray!10},major grid style={line width=.2pt,draw=gray!30}, xmajorgrids=true, ymajorgrids=true,  log ticks with fixed point, major tick length=0.05cm, minor tick length=0.0cm, legend style={at={(0.05,-0.115)}, legend image post style={scale=1}, column sep = 1pt, legend columns=2,draw=none}, scatter/classes={%
	non-optimal={mark=*,blue}, optimal={mark=triangle*,red}}]
  \nextgroupplot[%
        title= \footnotesize Iteration \footnotesize 1,
    ]
	\addplot[scatter,only marks,%
    		scatter src = explicit symbolic,
		    visualization depends on=\thisrow{size}\as\size,
            scatter/@pre marker code/.append style={
            /tikz/mark size=\size}]%
    	table[x=x, y=y, meta=label, col sep=comma] {figures/feats_1.txt};
  \nextgroupplot[%
        title=  {\footnotesize Iteration \footnotesize 5},
        yticklabels={,,}
    ]
	\addplot[scatter,only marks, 
    		scatter src = explicit symbolic,
		    visualization depends on=\thisrow{size}\as\size,
            scatter/@pre marker code/.append style={
            /tikz/mark size=\size}]%
    	table[x=x, y=y,meta=label,col sep=comma] {figures/feats_5.txt};
    \legend{ \footnotesize non-optimal, \footnotesize optimal}
    \end{groupplot} 
    \end{tikzpicture}
    \caption{The benefit of employing adaptive solution prediction on a TSP instance with 30 nodes. The axes are two statistical features described in Section \ref{subsection: features}, and each mark represents an edge in the TSP. The red marks denote the edges in the optimal (or shortest) tour, while the blue marks denote the non-optimal ones. The optimal and non-optimal edges are more separable after several iterations of ASP.} 
    \label{fig:feats}
\end{figure}

We provide an example to demonstrate the benefit of ASP in Figure~\ref{fig:feats}. The two sub-figures show the localization of edges for a TSP with respect to the two statistical features at different iterations of ASP. In the first iteration, it is difficult to separate the edges in the shortest tour (red marks) from those non-optimal edges (blue marks). In the subsequent iterations of ASP, edges are better characterized by the statistical features, because more accurate statistical information is extracted from better-quality solutions found by the heuristic search. After five iterations, the optimal edges have been shifted to the upper-right corner, while the non-optimal ones have moved to the lower-left corner. As a result, the ML model can better predict the optimal edges, and so is the optimal tour.

\subsection{Literature Review}
Recent years have witnessed many studies leveraging ML to develop heuristic methods for COPs, which can be broadly classified into reinforcement-learning-based (RL-based) and solution-prediction-based methods. The former trains an agent via RL to construct high-quality solutions directly~\citep{bello2016neural,khalil2017learning,nazari2018reinforcement,kool2018attention,deudon2018learning,xinSCZ21}, which is particularly effective for routing problems with solutions being an ordered sequence. Taking the TSP as an example, a RL-based approach typically forms feature representations for cities  (or node embedding) given their coordinates and then constructs a feasible tour by selecting the next city to visit based on the cities already visited. In particular, \citet{liang2021multi} propose to reconstruct node embedding for unvisited cities every several steps during the construction of a tour, of which the high-level idea is somewhat similar to ours.

Compared to RL-based methods, solution-prediction-based methods are more applicable to general COPs because of their flexibility in the use of downstream search methods. This type of method typically assumes easy instances of the problem under consideration have already been solved to optimality in advance. This makes it possible to train a ML model to learn the mapping from the features of a decision variable to its value in the optimal solution, supervised by solved problem instances with known optimal solutions. For instance, the ML model can be trained to classify whether an edge in the TSP is part of the optimal tour or not, based on the features (e.g., length) of that edge. To obtain accurate predictions, existing studies attempt to enrich the feature representation for decision variables. The most popular approach is learning feature representation using graph neural networks~(GNNs)~\citep{li2018combinatorial, selsam2018learning, zhang2020nlocalsat, joshi2019efficient,ding2020accelerating,selsam2019guiding}. GNNs characterize a decision variable by incorporating information from variables in the local neighborhood, which is obtained naturally for graph-based problems~\citep{li2018combinatorial} or can be defined manually according to constraint coefficient matrices for general COPs~\citep{ding2020accelerating}. Instead of learning feature representation using GNNs, \citet{sun2019using} define statistical measures as features that directly characterize a decision variable in terms of its prevalence in a set of randomly sampled solutions and their objective values.

Given an unseen COP instance, a ML model trained for this problem can be used to predict the values of its decision variables in the optimal solution. The prediction of the optimal solution, if sufficiently accurate, can be used to guide a search method to quickly find quality solutions. Several examples include fixing a proportion of decision variables to their predicted values to reduce the problem size for an optimization method~\citep{sun2019using,ding2020accelerating}; rounding the ML prediction to integer values to seed the initial position of a stochastic local search method ~\citep{zhang2020nlocalsat}; using the ML prediction to guide a tree search method to visit first the nodes (i.e., sub-problems) likely to contain optimal solutions~\citep{selsam2019guiding,li2018combinatorial}; or using the ML prediction to build a probabilistic model for sampling~\citep{shen2022enhancing}. It is apparent that both generic methods~\citep{ding2020accelerating,shen21learning} and specialized methods~\citep{selsam2018learning,sun2019using} can benefit from solution prediction. 

Apart from developing heuristic methods, ML has also been used to accelerate the branch-and-bound (B\&B) exact method~\citep{lodi2017learning} for deriving a certificate of optimality. In particular, recent studies~\citep{ quesnelWDS22,shen2022enhancing,abs-2201-02535} have applied solution-prediction-based methods to boost column generation~(CG)~\citep{lubbecke2005selected} via efficiently finding good solutions (i.e., columns with negative reduced costs) by solving a series of NP-hard pricing problems. The enhancement in CG can lead to substantial improvements in Branch-and-Price~(B\&P)~\citep{barnhart1998branch}, where CG is used to solve linear-programming relaxations during the computational process of B\&B. In addition, ML has been used to make critical decisions that may greatly influence the computational time of B\&B, such as approximating the strong branching decisions~\citep{lodi2017learning, gasseCFC019}, estimating the size of the B\&B tree to determine whether to perform a restart~\citep{hendel2021estimating}, deciding when to perform Dantzig-Wolfe decomposition~\citep{kruberLP17} or how to perform it~\citep{bassoCT20}, and others~\citep{vaclavikNSH18,morabitD021}. For a more comprehensive overview, please refer to \citep{bengio2020machine}.

\subsection{Contributions}
This paper explores the idea of leveraging feedback from the search of a specific test instance to better characterize its decision variables, in order to improve the prediction accuracy of an offline-trained ML model. Our main contributions include the following:
\begin{itemize}
    \item We propose a novel framework, adaptive solution prediction (ASP). ASP implements a feedback mechanism from heuristic search to the ML model using statistical features, which extract useful information from solutions found by the search method and inform the ML model as to which value of a decision variable is likely to be optimal. We design ASP to maintain a set of best solutions sampled according to ML predictions across iterations. The evolving sample set can provide more reliable statistical information on the optimal solution, which can in turn help the ML model make more accurate predictions to boost heuristic search. 
    
    \item We demonstrate the general applicability of our ASP approach using different ML models on several NP-hard COPs. Empirical results show that, by harnessing feedback from search, ASP can 1) significantly improve the prediction quality of a ML model, resulting in much better heuristic solutions, and 2) better scale to larger problem instances.
    
    \item We further demonstrate the efficacy of ASP in boosting CG and B\&P with an application to the graph coloring problem, where ASP is used as a heuristic-pricing method for solving the NP-hard pricing problem. Empirical results show that ASP improves the performance of CG and B\&P compared to existing pricing methods.
\end{itemize}

Our paper is organized as follows. Section~\ref{sec:asp-generic} presents the ASP framework for tackling general COPs. Section~\ref{sec:asp-specific} demonstrates how ASP can be applied to three tested COPs, followed by computational results in Section~\ref{sec:ret}. Section~\ref{sec:cgbp} applies ASP to enhance CG and the B\&P exact method. We conclude our paper in Section~\ref{sec:conclusion}.

%% file: 02-approach.tex
\section{Adaptive Solution Prediction Framework}
\label{sec:asp-generic}
We describe the general methodology of ASP for solving COPs that can be formulated as MIP problems with binary decision variables. This section is organized following the workflow of ASP illustrated in Figure~\ref{fig:asp}. We describe the training phase of ASP in Section~\ref{subsec:training}. For the testing phase, we first present a search method that makes use of the ML prediction for sampling good solutions in Section~\ref{subsection: ps}; then we describe the feedback mechanism from the search method to the ML model via statistical features in Section~\ref{subsection: features}. Combining the above components, we outline the ASP algorithm and discuss the design of ASP in Section~\ref{subsec:asp}. Lastly, we discuss the relationships between ASP and conventional metaheuristics in Section~\ref{subsec:MHcompare}.

\subsection{Model Training and Prediction}
\label{subsec:training}
Following previous work \citep{li2018combinatorial,lauri2019fine,sun2019using,ding2020accelerating,zhang2020nlocalsat,abbasi2020predicting}, we construct the training data from a set of solved problem instances of a COP, each of which is associated with an optimal solution. In the training data, a training example corresponds to a decision variable $x_i$ in a problem instance, which contains a feature vector $\bm{f}_i$ and a label $y_i$. The feature vector~$\bm{f}_i$ characterizes the properties of that variable and the label $y_i$ holds its value in the optimal solution associated with that problem instance. We employ two types of features including statistical features and problem-specific features. We describe the statistical features in Section~\ref{subsection: features} and the problem-specific ones when introducing the test problems in Section~\ref{sec:asp-specific}.

Using the training data, an off-the-shelf ML model can be trained to best separate the training examples with labels $y=1$ and those with labels $y=0$, based on their feature representations $\bm{f}$. Take the support vector machine with a linear kernel~(linear-SVM)~\citep{boser1992training} as an example. The linear-SVM aims to find the maximum margin hyperplane (parameterized by $\bm{w}$ and $b$) by solving the following quadratic programming problem with linear constraints:
\begin{align*}
	\min_{\bm{w},b, \bm{\xi}} & \quad \frac{1}{2}\bm{w}^T\bm{w} + r^+ \sum_{y_i=1}\xi^i + r^- \sum_{y_i=0}\xi^i, \\ 
	s.t. & \quad y_i\big(\bm{w}^T\bm{f}_i+b\big) \ge 1 - \xi^i, \quad i = 1, \cdots n_t,\\
	& \quad \xi^i \ge 0, \quad i = 1, \cdots n_t,
\end{align*}
where $r^+>0$ and $r^->0$ are respectively the penalty parameters for misclassifying positive ($y_i = 1$) and negative ($y_i = 0$) training examples; $\xi^i$ are the slack variables; and $n_t$ is the number of training examples.

In the testing phase, the offline-trained linear-SVM can be used to predict the optimal solution values of the decision variables for an unsolved problem instance. More specifically, for a decision variable $x_i$, we can compute the distance from the feature vector $\bm{f}_i$ associated with it to the optimal decision boundary of the linear-SVM $(\bm{w}^*,b^*)$ in the feature space:
\begin{equation}
d_i=\frac{{\bm{w}^*}^{T}\bm{f}_i+b^*}{||\bm{w}^*||}.  
\end{equation}
The signed distance $d_i$ indicates which value this decision variable is likely to take in the optimal solution. In an extreme case where $d_i$ approaches infinity, the linear-SVM is confident to predict the optimal solution value of this decision variable to be~$1$. On the other hand, if $d_i$ approaches negative infinity, the optimal solution value of this variable is highly likely to be $0$. We then use a logistic function (parameterized by $\beta_0$ and $\beta_1$), 
\begin{equation}
    p(d_i; \beta_0, \beta_1) = \frac{1}{1+e^{-(\beta_0 d_i + \beta_1)}},
\end{equation}
to normalize the distance value $d_i$ into a real number $p_i \in [0,1]$. For the unsolved instance, we can interpret the ML prediction $\bm{p} = [p_1, p_2, \cdots, p_n]$ of its optimal solution as the probabilities of the decision variables taking the value of $1$ in the optimal solution, 

We illustrate the above procedure using the maximum weight clique problem (MWCP) as an example. The MWCP aims at finding a clique on a graph such that the sum of weights associated with the vertices in the clique is maximized. The MWCP can be formulated as a MIP with binary decision variables indicating whether a vertex is used or not. We train a linear-SVM using data formed from a set of optimally solved MWCP instances. In the training data, a training example contains features (e.g., the degree of the vertex) representing the corresponding vertex in a MWCP instance and a label indicating whether the vertex is in the optimal clique to the MWCP instance. Given an unseen MWCP instance, the offline-trained linear-SVM can be used to predict whether a vertex is likely in the optimal clique, and the ML prediction is rescaled by a logistic model to produce probabilistic outputs, i.e., the probabilities of vertices used in the optimal clique.

\subsection{Probabilistic Sampling as a Search Method}
\label{subsection: ps}

\begin{figure}[th!]
    \centering
\scalebox{1}{\begin{algorithm}[H]
\LinesNumbered
\SetAlgoLined
\KwIn{$\mathcal{P}(\bm{x};\bm{c},\bm{A}, \bm{b})$: problem data;\\
 \hspace{3em} $\bm{p}$: the ML prediction of the optimal solution.}

   $\mathcal{S} \gets \varnothing$, and $\mathcal{C} \gets \{v_1, \cdots, v_n\}$\;
   
   \While{$\mathcal{C}$ is not empty}{
        $i \sim \frac{p_i}{\sum_{j \in \mathcal{C}} p_j}$, where $i \in \mathcal{C}$\; 
        $\mathcal{S} \gets \mathcal{S} \cup x_i$\;
        $\mathcal{C} \gets \mathcal{C} \setminus conflict(\mathcal{S})$\;
   }
    
    \Return{$\mathcal{S}$}
 \caption{Probabilistic Sampling.}
 \label{alg:ps}
\end{algorithm}}
\end{figure}

We use the ML prediction to form a probabilistic model to sample solutions, namely the probabilistic sampling (PS) method. Our PS method is outlined in Algorithm~\ref{alg:ps}. Let $\mathcal{S}$ denote a partial solution containing decision variables fixed to $1$ and let $\mathcal{C}$ denote the set of unfixed variables (Line~1). A solution is constructed by iteratively fixing a variable $x_i$ to $1$, and the variable is sampled from a probability distribution over all the variables in the candidate set~$\mathcal{C}$ (Line~3). Here, the probability of sampling a variable~$x_i \in \mathcal{C}$ is proportional to its value~$p_i$ predicted by the linear-SVM, and the candidate variables are those that can be added to the current partial solution~$\mathcal{S}$. Once a new variable is fixed to~$1$, the conflict variables are detected and removed from the candidate set~(Line 5). Note that this step is commonly referred to as domain propagation in the constraint-programming community, and generic techniques exist~\citep{AchterbergBKW08}. Finally, the process terminates and returns a feasible solution when all the variables are fixed, i.e., the set~$\mathcal{C}$ is empty (Line~2).

For the MWCP, the PS method constructs a clique by iteratively sampling promising vertices from the candidate set~$\mathcal{C}$ according to the ML prediction. To ensure the feasibility of the solution, the vertices in the candidate set must be neighbors to all the vertices already in the partial clique. We use an efficient algorithm based on the set representation of cliques~\citep{sun2019using} to maintain a valid candidate set.

\subsection{Refining ML Prediction via Statistical Features}
\label{subsection: features}

We improve the prediction quality of the offline-trained linear-SVM by harnessing the feedback from the PS method. Specifically, we make use of two statistical measures~\citep{sun2019using} as a part of the feature representation for a decision variable. Given sample solutions from the PS method, these statistical features can characterize a decision variable in terms of its prevalence in the sampled solutions and their objective values. Once the statistical features for decision variables are refined, an offline-trained ML model can be used to make a new prediction $\bm{p}$. For the MWCP, the statistical scores suggest whether a vertex frequently appears in high-quality cliques among samples, and inform the ML model that the vertex is more likely to be in the optimal clique.

The first statistical feature computes a ranking score,

\begin{equation}
    f_r(x_i) = \sum^{K}_{k=1} \frac{\bm{s}^{k}_{i}}{r^{k}},
\end{equation}

\noindent where $\bm{s}^{k}_{i}$ denotes the value of the decision variable~$x_i$ in the $k^{th}$ solution in a sample set $\mathbb{S}$; $r^k$ denotes the rank of the $k^{th}$ solution in terms of its objective value; and $K$ denotes the sample size~$ K=|\mathbb{S}|$. A variable~$x_i$ with a high ranking score indicates that 1) the variable appears more frequently in the sampled solutions, and 2) the ranks of those solutions are high. Since this ranking-based feature accumulates the scores from sampled solutions, it is sensitive to the size~$K$. Thus, we normalize it using the maximum ranking score among all variables in a problem instance.

The second statistical feature computes the Pearson correlation coefficient between the values of a variable $x_i$ and the objective values of the sampled solutions,
\begin{equation}
\label{eq:corr}
    f_c(x_i) = \frac{\sum^{K}_{k=1} (\bm{s}^{k}_{i} - \overline{\bm{s}}_{i}) ({o}^{k} - \overline{o})}{\sum^{K}_{k=1} \sqrt{(\bm{s}^{k}_{i} - \overline{\bm{s}}_{i})^2} \sqrt{\sum^{K}_{k=1} ({o}^{k} - \overline{o})^2}},
\end{equation}

\noindent where $\overline{\bm{s}}_{i}$ denotes the mean of $x_i$'s values across the sampled solutions; ${o}^{k}$ and $\overline{o}$ are the objective value of $k^{th}$ solution and the mean objective value of the sampled solutions in $\mathbb{S}$, respectively. We normalize the correlation-based measure according to the type of optimization problem. For a maximization problem, a variable highly positively correlated with the objective value is more likely to take the value of $1$ in high-quality solutions. Thus, this feature is normalized with respect to the maximum correlation score among variables. In a minimization problem, it is normalized by the minimum correlation score.

\subsection{Adaptive Solution Prediction Algorithm}
\label{subsec:asp}

\begin{figure}[th!]
    \centering

\scalebox{1}{\begin{algorithm}[H]
\LinesNumbered
\SetAlgoLined
\KwIn{$\mathcal{P}(\bm{x};\bm{c},\bm{A}, \bm{b})$: a problem instance; $\mathcal{ML}_{\bm{\theta}}$: an offline-trained ML model;\\
      \hspace{3em}$T$: number of iterations; $M$: sample size.\\}

    $\mathbb{S}^{(1)}$: uniformly and randomly sampled $M$ solutions\;
    \For{$t\gets1$ \KwTo $T$}{
      $\bm{F}^{(t)} \gets \bm{f}(\mathbb{S}^{(t)})$, where $\bm{F} = [\bm{f}^{(t)}_1, \cdots, \bm{f}^{(t)}_n]$\; 
      $\bm{p} \gets \mathcal{ML}(\bm{F}^{(t)}; \bm{\theta})$, where
      $p_i = P(x_i = 1)$, $i=1..n$\;
      \For{$m\gets1$ \KwTo $M$}{
      $\mathcal{S} \gets$ PS($\mathcal{P}$, $\bm{p}$)  \tcp*{Alg. \ref{alg:ps}}
      
      Let $\hat{\mathcal{S}}$ denote the worst solution in $(\mathbb{S}^{t})$\;
      \If{$\mathcal{S} \not\in \mathbb{S}^{t}$ and $obj(\hat{\mathcal{S}})$ worse than $obj(\mathcal{S})$}{
            replace $\hat{\mathcal{S}}$ with $\mathcal{S}$ in $\mathbb{S}^{t}$\;
       }
      
      }
    }
    
    \Return{$\bm{p}$}
 \caption{Adaptive Solution Prediction Algorithm}
 \label{alg:asp}
\end{algorithm}}
\end{figure}

Algorithm \ref{alg:asp} outlines our ASP method. In the beginning, a set of feasible solutions is generated uniformly and randomly (Line 1). Statistical features are then computed using the sampled solutions to form the feature representation of decision variables~$\bm{F}$~(Line 3). Based on $\bm{F}$, the offline-trained linear-SVM with the logistic model is used to predict a `probability' value for each decision variable (Line 4). The ML prediction~$\bm{p}$ is then used to guide the PS method to construct $M$ solutions, and the unique and better-quality new samples are used to replace the worst solutions in the sample set~$\mathbb{S}$~(Lines~5-9). The above process (Lines~3-9) is repeated for $T$ iterations, and the improved ML prediction is returned in the end (Line 10). 

We design ASP with the aim to \textit{constructing and maintaining a diverse set of high-quality solutions}. Firstly, continuous improvement of the \textit{quality} of solutions in the sample set is vital, as the values of decision variables in high-quality solutions are likely to be highly correlated to the values in the optimal solution. As the search method continuously supplies better-quality solutions, more reliable statistical features are obtained. As a result, the offline-trained linear-SVM can make a better prediction for the optimal solution which in turn better guides the search. To this end, we replace worse solutions in the sample set with better-quality solutions. Secondly, maintaining a \textit{diverse} sample set allows ASP to better explore the solution space and mitigate the possibility of premature convergence of the ML prediction to a local optimum. Therefore, we employ the PS method to construct a diverse set of solutions, benefiting from the randomness in sampling. Finally, we note that the PS method works well for ASP in improving the prediction quality of an offline-trained ML model, however, it may be inadequate for generating (near-)optimal solutions for certain COPs. For instance, \citet{joshi2019learning} have shown that on large TSP instances a sampling-based method is insufficient to produce high-quality solutions. In this case, we can use the improved ML prediction provided by ASP to boost a search method better suitable for that problem, such as specialized methods~\citep{selsam2019guiding,sun2019using}.

\subsection{Relationship to Existing Metaheuristics}\label{subsec:MHcompare}

In the vast literature of metaheuristics, many methods use a form of long-term memory during the solving process, such as evolutionary algorithms that keep a set of high-quality solutions explicitly or EDAs that build a probabilistic model from these solutions~\citep{hauschild2011introduction}. Our proposed ASP is mostly similar to EDAs, as they both alternate between randomly sampling the solution space of a specific test problem instance and updating the probability distribution to increase the likelihood of producing `good' solutions. In this perspective, the feedback mechanism in ASP can be analogous to the online-learning mechanisms of conventional EDAs. However, ASP also benefits from offline learning via ML. Specifically, ASP employs a ML model trained offline under the supervision of historically solved problem instances. We illustrate this benefit of ASP by comparing it with two representatives in the EDA class as follows.

ASP can be more efficient than some EDAs that require refitting the model during the problem-solving process. For instance, Bayesian optimization algorithms~\citep{el2019updating} leverage the Bayesian network to capture the relationship between decision variables. During the running of the algorithm, the Bayesian network is refitted at every iteration using the current set of sampled solutions. This refitting is very expensive and is one of the limitations of this type of algorithm. On the other hand, ASP leverages an offline-trained ML model and does not refit the model during the problem-solving process.  

ASP can be naturally seen as a generalization of simple EDAs~\citep{dorigo1996ant}, i.e., those using a simple probability model and only updating the parameters rather than its structure. For instance, ant colony optimization (ACO)~\citep{dorigo1996ant} constructs a probabilistic model where the probability of setting a variable $x_i$ to one is proportional to $\tau_i^\alpha\eta_i^\beta$. $\tau$, the pheromone value, is essentially a statistical measure (similar to the statistical measures employed by ASP) computed from sampled solutions; $\eta_i$ is a heuristic weight designed by human experts to guide the solution; $\alpha$ and $\beta$ are the parameters weighing these measures. While ACO relies on the two measures and combines them with a manually-defined function to form probabilistic models, our ASP can incorporate many more useful measures (or features) and learn a model automatically via ML.

%% file: 03-problems.tex
\section{Adapative Solution Prediction for Three NP-hard COPs}\label{sec:asp-specific}

In this section, we demonstrate how to apply ASP to our test COPs, respectively. For each COP, we present the problem formulation, problem-specific features, and the specialized algorithm for constructing feasible solutions in the PS method.

\subsection{Maximum Weight Clique Problem}
\label{subsec:mwcp}
The MWCP aims at finding a clique on a graph such that the sum of vertex weights of that clique is maximized. Let $G (\mathcal{V}, \mathcal{E}, \mathcal{W})$ be an un-directed weighted graph, consisting of a set of vertices $i \in \mathcal{V}$, a set of edges $(i, j)\in \mathcal{E}$, and a set of weights $w_i \in \mathcal{W}$ each of which corresponding to a vertex $i$. Let $\overline{\mathcal{E}} = \{(i,j)\,|\,\forall (i,j)\not\in \mathcal{E}\}$ denote the set of edges not in $\mathcal{E}$ and let the decision variable $x_i \in \{0,1\}$ indicate whether a vertex $i$ is a part of the clique. The MWCP can be defined as

\begin{align}
    \max_{\bm{x}} \; & \sum_{i \in \mathcal{V}} w_i x_i, \\
    s.t.\; & x_i + x_j \le 1;\; && (i,j) \in \overline{\mathcal{E}}, \\ 
          & x_i \in \{0,1\};\; && i=\{1,2,\cdots, |\mathcal{V}|\}.
\end{align} 

To apply ASP to the MWCP, we train a linear-SVM to predict whether a vertex $i$ is a part of the maximum weight clique. To initialize the statistical features, we generate random cliques using PS with uniform probabilities. In addition, we adopt several features to characterize a vertex~$i$, including the weight of the vertex $w_i$; the degree of the vertex, defined by the number of its neighbors $|\mathcal{N}_i|$; the upper bound of the clique containing the vertex $i$, computed by $w_i + \sum_{j \in \mathcal{N}_i} w_j$; the graph density, computed by $2|\mathcal{E}|/(|\mathcal{V}|(|\mathcal{V}|-1))$. 

\subsection{Travelling Salesman Problem}
Given a set of $n$ cities, the TSP aims to find a tour visiting each city exactly once such that the distance of the tour is minimized. Let the decision variable $x_{i,j}\in \{0,1\}$ denote whether the edge $(i,j)$ is in a feasible tour; let $d_{i,j}$ denote the distance of the edge $(i,j)$; let $Q$ denote a subset of $n$ cities. The TSP can be modeled using the Dantzig–Fulkerson–Johnson formulation~\citep{dantzig1954solution}:
 \begin{align}
    \min_{\bm{x}} \;& \sum_{i=1}^{n}\sum_{j=1}^{n}d_{i,j}x_{i,j}, \\
        s.t. \;& \sum_{i=1}^{n} x_{i,k} = \sum_{j=1}^{n} x_{k,j} = 1;\; && 1\leq k \leq n,  \\
              & \sum_{i\in Q}\sum_{j\neq i, j\in Q} x_{i,j} \leq |Q| - 1; && \forall Q \subsetneq \{1,\cdots, n\}, |Q| \geq 2, \\
              & x_{i,j} \in \{0,1\};\; && 1 \leq i,j \leq n,
\end{align}
\noindent where Constraints~(9) specify that each city must be visited once and Constraints~(10) prevent sub-tours, i.e., sub-tour elimination constraints.

To apply ASP to the TSP, we train a linear-SVM to predict whether an edge is in the shortest tour, supervised by a set of solved small TSP instances with known shortest tours (i.e., optimal solutions). To initialize the statistical features, we use feasible tours randomly generated by permuting the order of cities. We derive problem-specific features to the TSP from the length $d_{i,j}$ of an edge $(i,j)$, normalized by the mean or minimum length of edges connected to the city $i$ or $j$. Let $\overline{d}(i) = \frac{1}{n} \times \sum_{k=1}^{n} d_{ik}$ be the average length of edges connected to $i$.
Let $d^{max}(i)=\max_{k=1,\cdots,n} d_{i,k}$ and $d^{min}(i)=\min_{k=1,\cdots,n} d_{i,k}$ be the edge of maximum length and that of minimum length connected to $i$, respectively. The four normalized features are computed by
\begin{equation}
\frac{d_{i,j} - \overline{d}(i)}{d^{max}(i)-d^{min}(i)}, \;\;\;\;\;
\frac{d_{i,j} - d^{min}(i)}{d^{max}(i)-d^{min}(i)},
\;\;\;\;\;
\frac{d_{i,j} - \overline{d}(j)}{d^{max}(j)-d^{min}(j)}, \;\;\;\;\;
\frac{d_{i,j} - d^{min}(j)}{d^{max}(j)-d^{min}(j)}. 
\end{equation}

Given an unseen TSP instance, the trained linear-SVM is used to predict whether an edge is likely to be in the shortest tour, and the ML prediction is used to guide the PS method for constructing new tours. Specifically, PS starts with the city indexed by $i \mod n$ ($i$ being the index to the current sample) and iteratively samples the next city to visit, given the current visiting city $i$ and past visited cities. The likelihood of visiting the city $j \in \mathcal{C}$ is proportional to the ML prediction $\bm{p}_j$ of the edge $(i, j)$. After all the cities are visited, PS returns to the first-visited city to form a feasible tour. When updating the sample set~$\mathbb{S}$, we keep half of the initial random samples at the subsequent iterations. This maintains the diversity of the sample set~$\mathbb{S}$ and results in a more steady improvement of the ML prediction.

Since PS is not sufficient for the solving of large TSP instances~\citep{joshi2019learning}, we employ a greedy method to make use of the ML prediction for constructing high-quality tours, in addition to PS. At each iteration of ASP, the greedy method constructs tours starting from different cities in a similar way to PS, except that the next visiting city (given the visited ones) is selected greedily rather than sampled randomly according to the ML prediction. 

\subsection{Orienteering Problem}

The orienteering problem (OP) aims to selectively visit a set of locations under some distance (or time) budget through a closed tour passing the depot location, such that the total prize collected from the visited locations is maximized. Let $i$ be the index to a set of $n$ locations associated with prizes~$p_i$; let $T$ denote the budget limit; let $d_{i,j}$ denote the distance of the edge~$(i,j)$; let $x_{i,j} \in \{0,1\}$ denote a decision variable that indicates whether the edge~$(i,j)$ is a part of the tour; let $u_i$ denote the visiting order of a location $i$.  Given the depot location indexed by $1$, the OP can be formulated as 

 \begin{align}
    \max_{\bm{x}, \bm{u}} \;& \sum_{i=2}^{n}\sum_{j=2}^{n}p_{j}x_{i,j}, \\
        s.t. \;& \sum_{j=2}^{n} x_{k,j} = \sum_{i=2}^{n} x_{i,k} \leq 1;\; && 2 \leq k \leq n,  \\
              & \sum_{j=1}^{n} x_{1,j}= \sum_{i=1}^{n} x_{i,1}=1, \\
              & \sum_{i=2}^{n}\sum_{j=2}^{n} d_{i,j}x_{i,j} \leq T,\\
              & \sum_{i\in Q}\sum_{j\neq i, j\in Q} x_{i,j} \leq |Q| - 1; && \forall Q \subsetneq \{2,\cdots, n\}, |Q| \geq 2, \\
              & x_{i,j} \in \{0,1\};\; && 1 \leq i,j \leq n,
\end{align}
\noindent where Constraints~(14) ensure the connectivity of a tour; Constraints (15) specify that the tour must start from and return to the depot, respectively; Constraint~(16) is the capacity constraint; Constraints (17) are the sub-sour elimination constraints same as those in the TSP.

To apply ASP to the OP, we train a linear-SVM to predict whether an edge $(i,j)$ is in the most profitable tour. To initialize the statistical features, random solutions are generated using PS (Algorithm~\ref{alg:ps}) by setting a uniform probability (e.g., $\bm{p}=\bm{1}$) for all decision variables. We employ three problem-specific features. The first one is the distance of the edge normalized by the budget limit, computed by $\frac{d_{i,j}}{T}$. The remaining two features are based on the ratio between the prize gained by visiting a location and the distance of traveling to the location given the current location. The two features are normalized differently, respectively defined as 
\begin{equation}
    \frac{p_j / d_{i,j}}{\max_{k=1,\cdots,n} p_k / d_{i,k}},\;\;\;\frac{p_j / d_{i,j}}{\max_{k=1,\cdots,n} p_k / d_{k,j}}.
\end{equation}

Given an unseen OP, the trained linear-SVM can predict whether the edges are likely in the most profitable tour, and the ML prediction is used to guide PS. Similar to the TSP, PS starts from the depot location and iteratively samples the next location from a candidate set to visit. Given the current partial tour, the candidate locations are those 1) that are not visited before, and 2) that meet the condition: the total distance from the current position to the candidate vertex and then to the depot is no larger than the remaining budget. PS returns to the depot location to form a feasible solution when the candidate set is empty.

\subsection{ASP settings}
\label{subsec:asp-param}

\begin{table}[th!]
\centering
\caption{Problem-specific features.}
\label{tab:pfeat}
\resizebox{0.75\columnwidth}{!}{
\begin{tabular}{@{}ccc@{}}
\toprule
Problem & Decision variable & Problem-specific features\\
\midrule
TSP & edge & normalized edge distances\\ \midrule
OP & edge & normalized edge distance, normalized per-cost gains \\ \midrule
MWCP & vertex & graph density, vertex weight, vertex degree, vertex upper bound \\
\bottomrule
\end{tabular}}
\end{table}

This part discusses a general methodology for configuring ASP. Firstly, the choice of the ML model can have an impact on ASP, as complex ML models may provide more accurate predictions by capturing non-linear relationships between features, whilst simple ML models can be more efficient and robust in terms of generalization. We empirically set the linear-SVM as our default ML model for ASP, among several considered ML models including a classification tree and a feedforward neural network. A comparison between ASPs using different ML models will be shown in Section~\ref{subsec:exp-ml}.

We form the training data for each COP using $100$ optimally solved small problem instances. A common issue is that the number of positive training examples (e.g., optimal edges for a TSP) is much smaller than the number of negative ones. This can have a negative impact on the training of a ML model, as in our situation correctly classifying positive training examples is more critical than correctly classifying negative ones, also noted in~\citep{quesnelWDS22}. For the linear-SVM, we address this issue by raising the value of the penalty parameter for misclassifying positive examples~$r^{+}$ to the inverse of the ratio between the number of positive and negative examples in the training data, whereas the penalty parameter for misclassifying negative examples~$r^{-}$ is set to~$1$. Note that other parameters in the linear-SVM are set as default values in~\citep{chang2011libsvm}. For training other ML models, we replicate positive training examples multiple times such that the number of positive examples is about the same as the number of negative ones. After training a ML model, we can then tune the parameters~$\beta$ in the logistic model, which may significantly affect the efficiency of the downstream PS method. We use Bayesian Optimization (BO) for parameter tuning \citep{snoek2012practical, fernando2014bo}, which runs ASP with different sets of~$\beta$ values over~$30$ training instances and identifies the $\beta$ parameters that yield the best solution quality.

ASP has several other parameters, as shown in Algorithm~\ref{alg:asp}. The iteration number ($T$) need not be determined in advance, and ASP can be stopped when no better solution can be found in successive iterations for the current problem instance. In our experiments, we terminate the algorithm when reaching a certain time limit that is set to the computational time of a compared method. We suggest setting the sample size ($M$) proportional to the problem size. Empirically, we set the sample size to  $n$ for the MWCP, $20n$ for the TSP, and $50n$ for the OP, where $n$ denotes the number of vertices (or cities) in a graph. The performance of ASP is fairly robust with respect to the sample size. However, we note that too large a sample size (e.g., $100n$ for the TSP) can slow down the progress of ASP, because a ML model can make fewer predictions under a certain time budget. On the other hand, a very small sample size (e.g., $n$ for the TSP) may result in early convergence of ASP because of the insufficient exploration of the search space.

%% file: 04-result.tex
\section{Computational Results and Analysis}
\label{sec:ret}

In this section, we demonstrate that 1) ASP can effectively improve the ML prediction (by a linear-SVM) so as to guide the search to find better-quality solutions (Section~\ref{subsec:asp-exp}); 2) compared with existing ML-based methods, ASP benefits from its online-learning (or feedback) mechanism and can better generalize to unseen instances of different characteristics and sizes than the training ones; 3) compared with conventional heuristic methods, ASP benefits from offline learning via ML and can find better-quality solutions more efficiently (Section~\ref{subsec:exp-trad}). Lastly, we show that ASP, as a generic framework, can work effectively with several classes of ML models on the MWCP (Section~\ref{subsec:exp-other}).

We generate problem instances in different sizes following related studies. TSP instances are generated by sampling a set of locations uniformly randomly in a unit square~\citep{bello2016neural,joshi2019efficient, kool2018attention, kobeaga2018efficient}. For the OP, we generate a set of locations the same way as the TSP; the distance budget~$T$ is uniformly sampled from $[\frac{d}{2}-1, \frac{d}{2}+1]$ to maximize the problem difficulty~\citep{fischetti1998solving}, where $d$ is the average objective value of optimal solutions for solving TSPs of a particular size; the prize attached to a location can be generated by a certain scheme chosen from three schemes, constant, uniform, and distance~\citep{fischetti1998solving}. For the MWCP, we generate graphs with a certain size and density using the well-known Barabási-Albert (BA) method~\citep{barabasi1999emergence} and the Erdős-Rényi (ER) method~\citep{erdHos1960evolution}, and we assign weights to graph vertices according to~\citep{jiang2018two,cai2016fast}.

We compare different methods using primal gap, i.e., the normalized difference between the best-found solution and the optimal solution for a problem instance, defined as
\begin{equation}
\text{primal gap} \coloneqq \left|\frac{obj^{best\_found}}{obj^{optimal}} - 1\right| \times 100\%.
\end{equation} 
In addition, we use average precision~\citep{zhu2004recall} to assess the quality of ML predictions for solution-prediction-based methods following previous work~\citep{ding2020accelerating,shen21learning}. Average precision is computed as the area under a curve that measures the trade-off between precision and recall at different decision thresholds. We note that both measures (i.e., primal gap and the average precision) are computed using the optimal solution of the test problem instance. When the optimal solution is unknown, we report the objective value of the best-found solution directly.

\subsection{Adaptive Solution Prediction}
\label{subsec:asp-exp}

\begin{figure}[ht!]
    \centering
	\begin{tikzpicture}
	\begin{groupplot}[group style = {group size = 3 by 1, horizontal sep = 40pt}, height=0.35\textwidth, width=0.35\textwidth, grid style={line width=.1pt, draw=gray!10},major grid style={line width=.2pt,draw=gray!30}, xlabel =  \footnotesize $k^{th}$ iteration of ASP, axis y line*=left, ylabel =  \footnotesize Average precision of the ML prediction, y label style={at={(axis description cs:-0.12,.5)},anchor=south},x label style={at={(axis description cs:0.5,-0.12)},anchor=north}, xtick = {1,5,10,15,20}, xmajorgrids=true, ymajorgrids=true,  major tick length=0.05cm, minor tick length=0.0cm, legend style={at={(0.6,0.1)},anchor=south, column sep = 1pt,font=\scriptsize, legend columns = 1,draw=none}, title style={yshift=-1.5ex,}]
   \nextgroupplot[%
    title=\footnotesize {MWCP (ER graph, $1000$ nodes, density $0.18$)},
    ]
	\addplot[color=red, mark=square*, mark size = 1, line width=0.30mm] table [x=x, y=ER_0.18, col sep=comma] {figures/mwc_ap.txt};
   \nextgroupplot[%
    title=\footnotesize {TSP ($100$ cities)}, 
    ]
	\addplot[color=red, mark=square*, mark size = 1, line width=0.30mm] table [x=x, y=y100, col sep=comma] {figures/tsp_ap.txt};
	\nextgroupplot[%
    title=\footnotesize OP ($50$ locations and `uniform' price),
    ]
	\addplot[color=red, mark=square*, mark size = 1, line width=0.30mm] table [x=x, y=y50unif, col sep=comma] {figures/op_ap.txt}; 
    \end{groupplot} 
	\end{tikzpicture}
\caption{The average precision of the ML prediction at every iteration of ASP for three tested COPs. For each COP, the mean statistics across $100$ test instances are reported, and information for the test instances is shown in brackets. In ASP, the linear-SVM and the logistic model are trained using the solved problem instances of the same size as the test ones. }
\label{fig:asp-ap}
\end{figure}
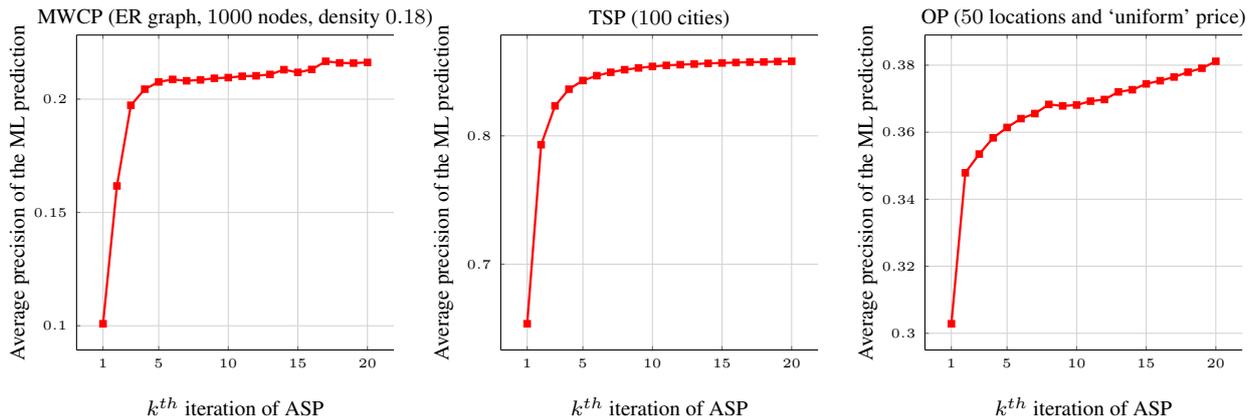

\def\paragraph#1{\noindent\emph{#1} }
\paragraph{Can ASP improve the quality of ML predictions?} Figure~\ref{fig:asp-ap} shows the average precision of the ML predictions over iterations of ASP. For each COP, ASP is trained using a certain size of problem instances and is tested on $100$ instances of the same size. As can be seen, on all three tested problems ASP can effectively improve the quality of the ML predictions, particularly in the first several iterations. Note that the improvement in prediction quality can be attributed to the feedback mechanism of ASP that recomputes statistical features over better-quality solutions. 

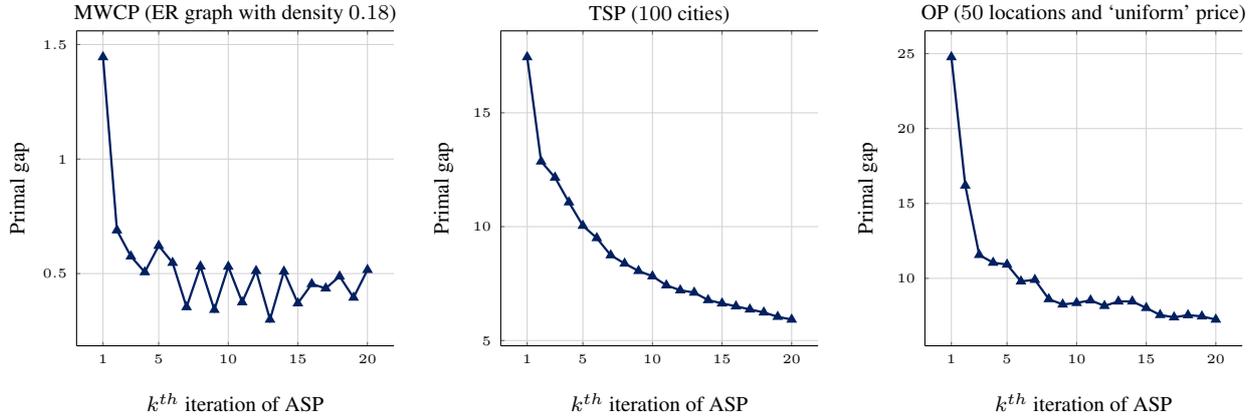
\begin{figure}[ht!]
    \centering
	\begin{tikzpicture}
	\begin{groupplot}[group style = {group size = 3 by 1, horizontal sep = 40pt}, height=0.35\textwidth, width=0.35\textwidth, grid style={line width=.1pt, draw=gray!10},major grid style={line width=.2pt,draw=gray!30}, xlabel =  \footnotesize $k^{th}$ iteration of ASP, axis y line*=left, ylabel =  \footnotesize Primal gap, y label style={at={(axis description cs:-0.12,.5)},anchor=south},x label style={at={(axis description cs:0.5,-0.12)},anchor=north}, xtick = {1,5,10,15,20}, xmajorgrids=true, ymajorgrids=true,  major tick length=0.05cm, minor tick length=0.0cm, legend style={at={(0.8,0.8)},anchor=south, column sep = 1pt,font=\scriptsize, legend columns = 1,draw=none}, title style={yshift=-1.5ex,}]
   \nextgroupplot[%
    title=\footnotesize MWCP (ER graph with density $0.18$),
    ]
	\addplot[color=marine, mark=triangle*, mark size = 1.5, line width=0.30mm] table [x=x, y=ER_0.18, col sep=comma] {figures/mwc_gap.txt};
   \nextgroupplot[%
    title=\footnotesize TSP ($100$ cities),
    ]
	\addplot[color=marine, mark=triangle*, mark size = 1.5, line width=0.30mm] table [x=x, y=y100, col sep=comma] {figures/tsp_greedy_obj.txt};

	\nextgroupplot[%
    title=\footnotesize OP ($50$ locations and `uniform' price),
    ]
    \addplot[color=marine, mark=triangle*, mark size = 1.5, line width=0.30mm] table [x=x, y=y50unif, col sep=comma] {figures/op_greedy_obj_opt.txt};
    \end{groupplot} 
	\end{tikzpicture}
\caption{The primal gap of the best solution at every iteration of ASP for three tested COPs. For each COP, the mean statistics across $100$ test instances are reported, and the size of the test problem instances is shown in brackets. In ASP, the linear-SVM and the logistic model are trained using the solved problem instances of the same size as the test ones. Note that the primal-gap curve is not guaranteed to monotonically decrease, as we report the best solution constructed using the ML prediction at a particular iteration rather than the best-found solution so far.}
\label{fig:asp-gap}
\end{figure}

\paragraph{Can improved ML predictions guide the search to find better solutions?}
Figure~\ref{fig:asp-gap} shows the primal gap of the best solution at each iteration of ASP. We can observe that the search method can construct much better solutions as a result of making use of the iteratively improved ML prediction. Furthermore, the improvements in the solution quality are particularly large in the first several iterations. This observation is consistent with the fast improvement of the prediction quality as shown in Figure~\ref{fig:asp-ap}. These results show that the effectiveness of the ML prediction and the quality of the best-found solution are closely related, and improving the quality of ML prediction is critical in boosting a search method. 

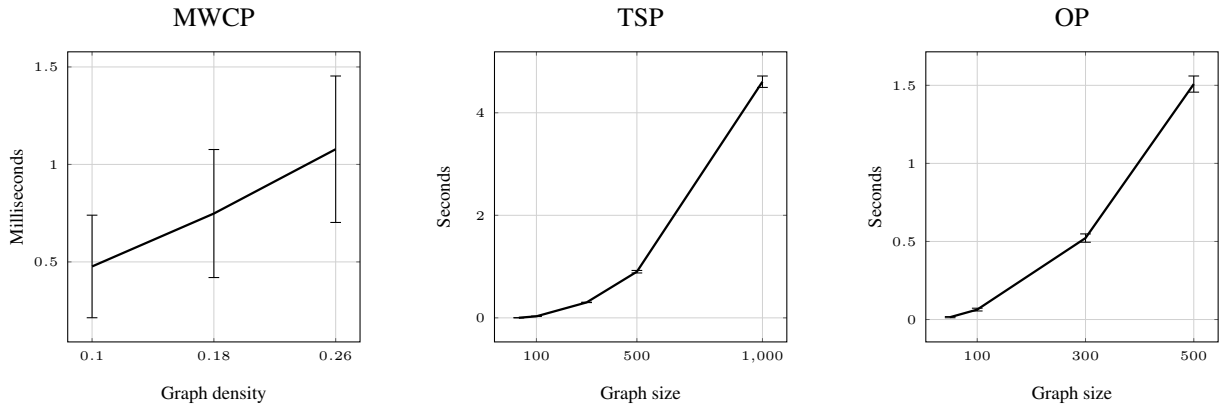
\begin{figure}[ht!]
    \centering
	\begin{tikzpicture}
	\begin{axis} [legend style={nodes={scale=0.5, transform shape}}, title=  MWCP, box plot width=0.0em,  height=0.33\textwidth,width=0.33\textwidth, grid style={line width=.1pt, draw=gray!10},major grid style={line width=.2pt,draw=gray!30}, xlabel = \scriptsize Graph density, ylabel = \scriptsize Milliseconds, y label style={at={(axis description cs:-0.12,.5)},anchor=south}, xtick = {0.1, 0.18, 0.26},x label style={at={(axis description cs:0.5,-0.12)},anchor=north}, xmajorgrids=true, ymajorgrids=true,  major tick length=0.05cm, minor tick length=0.0cm, legend style={at={(0.47,0.22)},anchor=west,font=,draw=none},
	mark repeat=1,mark phase=0]
	\addplot[color=black, line width=0.30mm, error bars/.cd, y dir=both, y explicit] table [x=name, y=mean, y error=std, col sep=comma] {figures/mwc_time.txt}; 
	\end{axis}
	\end{tikzpicture}
	\hspace{2em}
	\begin{tikzpicture}
	\begin{axis} [box plot width=0.0em, title=  TSP,  height=0.33\textwidth,width=0.33\textwidth, grid style={line width=.1pt, draw=gray!10},major grid style={line width=.2pt,draw=gray!30}, xlabel = \scriptsize Graph size, ylabel = \scriptsize Seconds, y label style={at={(axis description cs:-0.12,.5)},anchor=south}, xtick = {100,500, 1000},x label style={at={(axis description cs:0.5,-0.12)},anchor=north}, xmajorgrids=true, ymajorgrids=true,  major tick length=0.05cm, minor tick length=0.0cm, legend style={at={(0.6,0.24)},anchor=west,font=,draw=none},
	mark repeat=1,mark phase=0]
	\addplot[color=black, line width=0.30mm, error bars/.cd, y dir=both, y explicit] table [x=name, y=mean, y error=std, col sep=comma] {figures/tsp_time.txt}; 
	\end{axis}
	\end{tikzpicture}
	\hspace{2em}
	\begin{tikzpicture}
	\begin{axis} [box plot width=0.0em,  height=0.33\textwidth,width=0.33\textwidth, title=  OP, grid style={line width=.1pt, draw=gray!10},major grid style={line width=.2pt,draw=gray!30}, xlabel = \scriptsize Graph size, ylabel = \scriptsize Seconds, y label style={at={(axis description cs:-0.12,.5)},anchor=south}, xtick = {100,300,500},x label style={at={(axis description cs:0.5,-0.12)},anchor=north}, xmajorgrids=true, ymajorgrids=true,  major tick length=0.05cm, minor tick length=0.0cm, legend style={at={(0.6,0.24)},anchor=west,font=,draw=none},
	mark repeat=1,mark phase=0]
	\addplot[color=black, line width=0.30mm, error bars/.cd, y dir=both, y explicit] table [x=name, y=mean, y error=std, col sep=comma] {figures/op_time.txt};
	\end{axis}
	\end{tikzpicture}
\caption{Time spent per iteration when problem size grows.}
\label{fig:time}
\end{figure}

\paragraph{Computational time for ASP.} Figure~\ref{fig:time} shows the computational time spent per ASP iteration when the problem size grows. On the MWCP, ASP only spends a few milliseconds per iteration and its computational time increases slowly with respect to the number of edges in a graph. For the TSP and OP, the computational time per ASP iteration is in seconds and increases in polynomial time with respect to the number of cities. These observations can be explained by analyzing the PS method, which takes the majority of the computational time of ASP in practice. The time complexity of PS is determined by the sample size, the solution length, and the time complexity for resolving conflict variables. For the MWCP, its sample size and the length of a solution are much smaller than those for the TSP and OP. When problem size increases, the solution length for the MWCP grows much slower than that for the TSP and OP.

\subsection{Comparison to ML-based methods}
\label{subsec:exp-ml}

In this part, we compare ASP with existing ML-based methods, and show that the online-learning mechanism of ASP allows it to better generalize to unseen test problem instances. A method is trained using problem instances of a certain characteristic (e.g., problem size) and tested on instances of the same problem but with different characteristics. It is worth noting that scaling ML-based methods to larger problem instances has been a known challenge for the TSP~\citep{bengio2020machine}.

We consider the following ML-based methods as baselines: 1) single-shot solution prediction (SSSP), a special case of ASP that predicts the optimal solution only once at the initial iteration $T=1$. Since SSSP does not need to recompute statistical features and make predictions multiple times, it can construct more solutions than ASP under a certain cutoff time. 2) SSSP-GNN, which refers to a variant of SSSP that use the popular Graph Neural Network~(GNN) to predict optimal solutions. We consider the SSSP-GNN\footnote{Code: \url{https://github.com/chaitjo/graph-convnet-tsp}.} proposed by \citet{joshi2019learning} for the TSP and the one\footnote{Code: \url{https://github.com/intel-isl/NPHard}.} proposed by \citet{li2018combinatorial} for the MWCP, respectively. These methods are different in the model structure and the use of (problem-specific) search methods. We note that \citet{li2018combinatorial} initially tested their approach on a special case of the MWCP with the uniform vertex weight (i.e., $w=1$), hence they use a one-hot feature vector $\bm{x}=\bm{1}$ to represent a vertex. In our case, vertices can have different weights, therefore we set the feature vector for a vertex $i$ according to its weight, i.e., $\bm{x_i}=\bm{w_i}$. 3) RL-AM\footnote{Code: \url{https://github.com/wouterkool/attention-learn-to-route}.}, a RL-based approach designed for routing problems~\citep{kool2018attention}. RL-AM employs an attention-based neural model (AM) as the policy, which is trained end-to-end via RL to construct high-quality solutions directly.  

For testing, we parallelize the sampling process of our ASP and SSSP on $8$ CPUs (Intel(R) Core(TM) i5-8300H CPU @ 2.30GHz). For other methods with deep neural networks, we equip them with a Tesla P100 GPU, and their computational time can be substantially reduced as compared to that when using $8$ CPUs.

\begin{table*}[ht!]
\centering
\caption{The results of the ML-based methods on the MWCP. For a type of graph, the training set is formed from solved instances with $1000$ vertices and a density value of $0.18$, and the test instances are formed on the same type of graph but of different density values. A method is set to generate a total number of $10n$ solutions. The primal gap of the best solution is reported, averaged over $100$ test instances. The best results are highlighted in bold (statistically significant based on a $t$-test, with the significance level set to $0.05$).}
\label{tab:mwc}
\resizebox{0.95\textwidth}{!}{
\begin{tabular}{@{}cccccccccc@{}}
\toprule
\multirow{2}{*}{Graph type} & \multirow{2}{*}{Density} & \multicolumn{3}{c}{Average precision} & \multicolumn{3}{c}{Primal gap} & \multicolumn{2}{c}{Time} \\
& & ASP & SSSP & SSSP-GNN & ASP & SSSP & SSSP-GNN & ASP/SSSP & SSSP-GNN \\
\cmidrule(lr){1-2}\cmidrule(lr){3-5}\cmidrule(lr){6-8}\cmidrule(lr){9-10}
\multirow{3}{*}{Erdős-Rényi}& $0.10$ & \textbf{0.17} & 0.08 & 0.09 & 0.0 & 0.04 & 1.19 & 0.1 & 0.1\\
& $0.18$ & \textbf{0.23} & 0.14 & 0.10 & 0.0 & 0.0 & 3.74 & 0.1 & 0.1 \\
& $0.26$  & \textbf{0.21} & 0.12 & 0.09 & \textbf{0.17} & 0.23 & 6.71 & 0.1 & 0.2\\ 
\cmidrule(lr){1-2}\cmidrule(lr){3-5}\cmidrule(lr){6-8}\cmidrule(lr){9-10}
\multirow{3}{*}{Barabási-Albert}& $0.10$ & \textbf{0.95} & 0.84 & 0.49 & 0.0 & 0.0 & 3.78 & 0.1 & 0.2\\
& $0.18$  & \textbf{0.98} & 0.93 & 0.88 & 0.0 & 0.03 & 0.14 & 0.4 & 0.5\\
& $0.26$ & \textbf{0.99} & 0.97 & 0.53 & 0.0 & 0.06 & 8.82 & 0.2 & 0.3 \\ 
\bottomrule
\end{tabular}}
\end{table*}

Table~\ref{tab:mwc} shows the results on the MWCP. It can be seen that ASP produces the best-quality ML prediction among the three solution-prediction-based methods. Further, we observe the performances of both ASP and SSSP are robust, i.e., their performance does not significantly degrade when testing on the graph with a different density to the training ones. Being consistent with the average precision, the best-found solution by ASP has better quality than that of other methods.

\begin{table*}[ht!]
\centering
\caption{The results of the ML-based methods on the TSP. All methods are trained using TSP instances of $100$ cities, and are tested on TSP instances of different sizes. The cutoff time for ASP and SSSP is set to the minimum computational time between SSSP-GNN and RL-AM under their default parameter settings. The primal gap of the best solution is reported, averaged over $100$ test instances. The best results are highlighted in bold. }
\label{tab:tsp}
\resizebox{0.95\textwidth}{!}{
\begin{tabular}{@{}cccccccccccc@{}}
\toprule
\multirow{2}{*}{Graph size}& \multicolumn{3}{c}{Average precision}& \multicolumn{4}{c}{Primal gap} & \multicolumn{3}{c}{Time} \\
& ASP & SSSP & SSSP-GNN & ASP & SSSP & SSSP-GNN & RL-AM & ASP/SSSP & SSSP-GNN & RL-AM \\
\cmidrule(lr){1-1} \cmidrule(lr){2-4}\cmidrule(lr){5-8}\cmidrule(lr){9-11}
$30$ & \textbf{0.84} & 0.75& 0.71 & \textbf{0.6} & 1.0 & 33.6 & 2.2 & 0.1 & 0.1 & 0.1 \\
$100$ & 0.85 & 0.65& \textbf{0.98} & 4.1 & 11.9 & \textbf{2.1}  & \textbf{2.3} & 0.4 & 0.4 & 0.7 \\ 
$300$ & \textbf{0.85} & 0.52 & 0.42 & 12.0 & 32.9 & 71.5 & \textbf{10.4} & 3.8 & 3.8 &5.8 \\
$500$ & \textbf{0.84} & 0.46 & 0.34 & \textbf{14.1} & 41.5& 79.1 & 17.6 & 12.4 & 12.4 & 16.7 \\
$1000$ & \textbf{0.83} & 0.36 & 0.23 & \textbf{16.8} & 55.8 & 117.7 & 29.0 &48.6 &48.6 & 75.7 \\
\bottomrule
\end{tabular}}
\end{table*}

Table~\ref{tab:tsp} presents the results on the TSP, showing that ASP scales well across TSP instances of different sizes as evidenced by the average precision, while SSSP performs reasonably well only for small-scale TSP instances with 30 cities. Consequently, the quality of the best solution found by ASP is much more close to the optimal one than that found by SSSP. It is interesting to note that for SSSP-GNN and RL-AM, their performances are competitive only when the size of the training instances (i.e., $100$ cities) and testing instances are similar. However, their performances deteriorate quickly when testing on larger TSP instances.

\begin{table*}[ht!]
\centering
\caption{The results of the ML-based methods on the OP. All methods are trained using OP instances of $50$ locations and the uniform prize and are tested on OP instances of different sizes and different prize types. The cutoff time for ASP and SSSP is set to the average computational time of RL-AM under its default parameter setting. The OP is a maximization problem, and the mean objective value of best solutions averaged over $100$ test instances is reported. The best results are highlighted in bold.}
\label{tab:op}
\resizebox{0.95\textwidth}{!}{
\begin{tabular}{@{}cccccccccccc@{}}
\toprule
\multirow{2}{*}{Graph size} & \multicolumn{3}{c}{Prize type: `uniform'} & \multicolumn{3}{c}{Prize type: `constant'} & \multicolumn{3}{c}{Prize type: `distance'} & \multicolumn{2}{c}{Time}\\
& ASP& SSSP& RL-AM & ASP & SSSP & RL-AM & ASP & SSSP &RL-AM& ASP/SSSP & RL-AM\\
\cmidrule(lr){1-1} \cmidrule(lr){2-4}\cmidrule(lr){5-7}\cmidrule(lr){8-10}\cmidrule(lr){11-12}

$50$  & 14.9 & 14.8 & 15.0 & 27.1 & 26.9 & 26.5 & 14.9 & 14.7 & 14.9 & 0.2 & 0.2 \\
$100$ & 30.8 & 29.2 & 31.2 & \textbf{56.4} & 53.8 & 54.9 & 30.8 & 29.2 & 30.7 & 0.9 & 0.9\\
$300$ & 79.4 & 70.7 & \textbf{81.6} & \textbf{147.3} & 133.9 & 143.3 & \textbf{83.1} & 72.5 & 79.0 & 6.5 & 6.5\\
$500$ & \textbf{125.0} & 112.8 & 111.9 & \textbf{232.7} & 218.0 & 201.2 & \textbf{131.5} & 105.3 & 100.9 & 14.4 & 14.4\\

\bottomrule
\end{tabular}}
\end{table*}

Table~\ref{tab:op} shows the results on the OP. Overall, ASP significantly outperforms other methods in most cases, though RL-AM wins with statistical significance on just one occasion. Similar to the observations on the TSP, ASP significantly improves on SSSP, especially for large OP instances. Considering that the underlining distribution of real-world problems is usually unknown, these results demonstrate that ASP is more generally applicable to real-world scenarios.

Overall, several common observations can be made across the three COPs: 1) ASP can find much better solutions than its special case, SSSP. Given SSSP constructs more solutions than ASP, the observation shows that improving the prediction accuracy can significantly improve the efficiency of the search (i.e., PS). 2) Among all ML-based methods, ASP achieves much better scalability and/or generalization performance. Both observations should be attributed to the online-learning mechanism of ASP. Specifically, ASP improves the prediction of an offline-trained ML model by providing it with more accurate statistical information for characterizing decision variables. The information is extracted from better-quality solutions found \textit{during the search process of a specific test problem instance}. On the other hand, existing ML-based methods (e.g., SSSPs or RL-based methods) cannot leverage such instance-specific information.

\subsection{Comparison to traditional methods}
\label{subsec:exp-trad}


In this part, we compare the proposed ASP with traditional methods (not using ML) and show that ASP with the offline-learning mechanism (i.e., the use of ML) significantly outperforms two generic methods on all test COPs. Moreover, ASP achieves very competitive performances as compared with some specialized methods. For generic methods, we opt for ant colony system~(ACS)~\citep{dorigo1996ant} as the primary baseline, given the close relationship between ASP and EDAs~(Section~\ref{subsec:MHcompare}). ACS has been extensively studied in combinatorial optimization, and we adopt ACS variants in ~\citep{xu2007improved, ke2008ants, dorigo1997ant} for our test problems. In addition, we include Gurobi \citep{gurobi2018gurobi}, a commercial MIP solver with a bag of primal heuristics. As our focus is on solution quality, we instruct Gurobi to spend $90\%$ of the computational time on executing these primal heuristics by setting the `Heuristics' parameter to $0.9$. When using Gurobi to solve the TSP or the OP, the subtour elimination constraints (e.g., Constraints (7) for TSP) are implemented as lazy constraints, i.e., they are only added to a problem formulation when violated by the new infeasible solution.

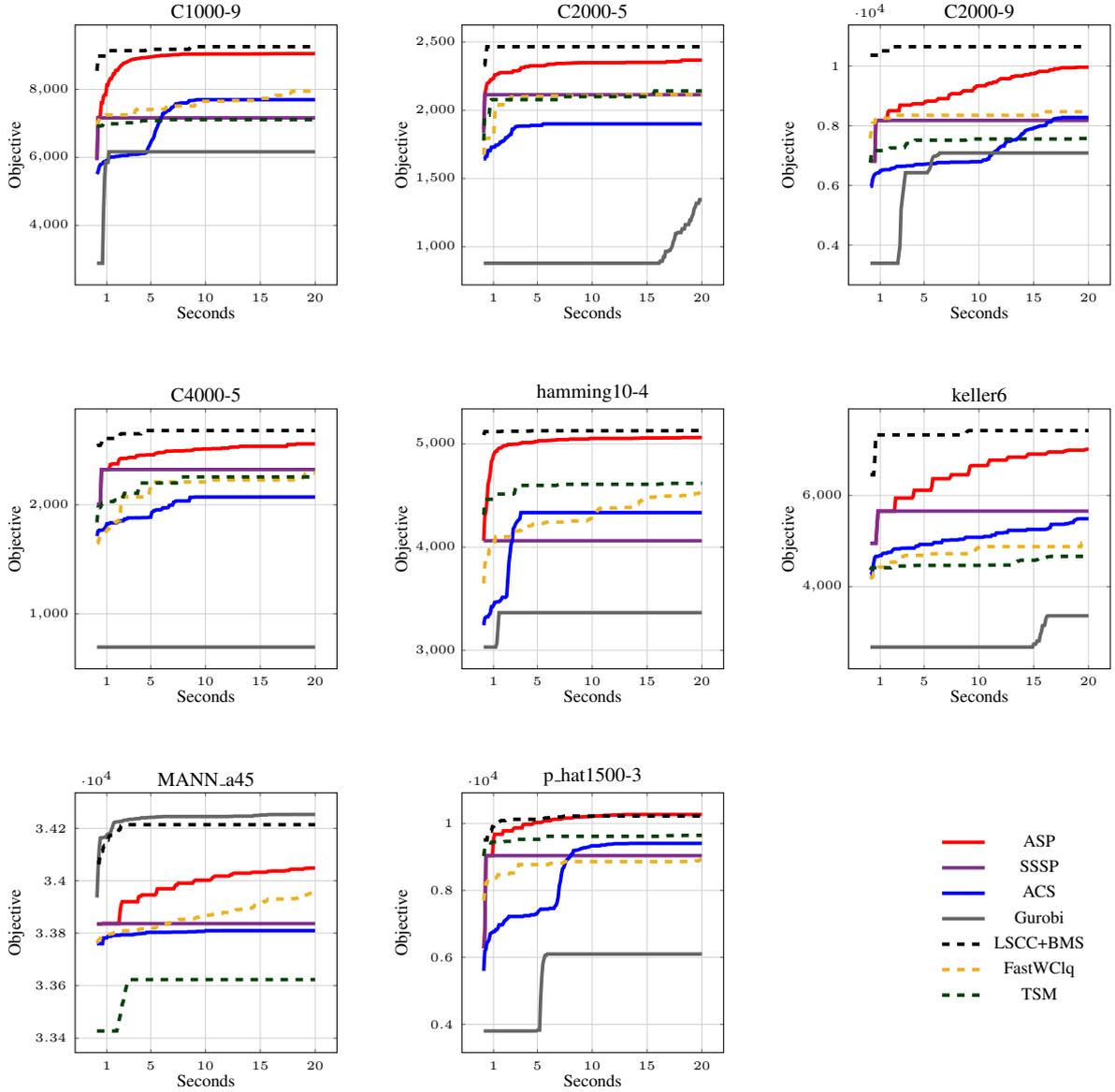
\begin{figure*}[ht!]
    \centering
	\begin{tikzpicture}
	\begin{groupplot}[group style = {group size = 3 by 3, horizontal sep = 50pt, vertical sep = 50pt}, height=0.32\textwidth, width=0.32\textwidth, grid style={line width=.1pt, draw=gray!10},major grid style={line width=.2pt,draw=gray!30}, xlabel = \scriptsize Seconds, ylabel = \scriptsize Objective, y label style={at={(axis description cs:-0.16,.5)},anchor=south}, xtick = {1,5,10,15,20},x label style={at={(axis description cs:0.5,-0.05)},anchor=north}, xmajorgrids=true, ymajorgrids=true,  major tick length=0.05cm, minor tick length=0.0cm, legend style={column sep = 1pt,font=\scriptsize, legend columns = 1,draw=none},title style={yshift=-1.5ex}]
   \nextgroupplot[%
   legend to name=grouplegend1,
    title= \footnotesize C1000-9,
    ]
	\addplot[color=red, line width=0.50mm] table [x=x, y=y, col sep=comma] {dimacs-new/C1000-9/mssp-0.txt}; \addlegendentry{ASP}
	\addplot[color=purple, line width=0.50mm] table [x=x, y=y, col sep=comma] {dimacs-new/C1000-9/sssp-0.txt}; \addlegendentry{SSSP}
	\addplot[color=blue, line width=0.50mm] table [x=x, y=y, col sep=comma] {dimacs-new/C1000-9/aco.txt}; \addlegendentry{ACS}
	\addplot[color=gray, line width=0.50mm] table [x=x, y=y, col sep=comma] {dimacs-new/C1000-9/gurobi.txt}; \addlegendentry{Gurobi}
	\addplot[color=black, line width=0.50mm, dashed] table [x=x, y=y, col sep=comma] {dimacs-new/C1000-9/lscc.txt};  \addlegendentry{LSCC+BMS}
	\addplot[color=orange, line width=0.50mm, dashed] table [x=x, y=y, col sep=comma] {dimacs-new/C1000-9/fastwclq.txt};  \addlegendentry{FastWClq}
	\addplot[color=green, line width=0.50mm, dashed] table [x=x, y=y, col sep=comma] {dimacs-new/C1000-9/tsm.txt};\addlegendentry{TSM}
   \nextgroupplot[%
    title= \footnotesize C2000-5,
    ]
	\addplot[color=red, line width=0.50mm] table [x=x, y=y, col sep=comma] {dimacs-new/C2000-5/mssp-0.txt}; 
	\addplot[color=purple, line width=0.50mm] table [x=x, y=y, col sep=comma] {dimacs-new/C2000-5/sssp-0.txt}; 
	\addplot[color=blue, line width=0.50mm] table [x=x, y=y, col sep=comma] {dimacs-new/C2000-5/aco.txt}; 
	\addplot[color=gray, line width=0.50mm] table [x=x, y=y, col sep=comma] {dimacs-new/C2000-5/gurobi.txt}; 
	\addplot[color=black, line width=0.50mm, dashed] table [x=x, y=y, col sep=comma, dashed] {dimacs-new/C2000-5/lscc.txt};  
	\addplot[color=orange, line width=0.50mm, dashed] table [x=x, y=y, col sep=comma] {dimacs-new/C2000-5/fastwclq.txt};  
	\addplot[color=green, line width=0.50mm, dashed] table [x=x, y=y, col sep=comma] {dimacs-new/C2000-5/tsm.txt};
   \nextgroupplot[%
    title= \footnotesize C2000-9,
    ]
	\addplot[color=red, line width=0.50mm] table [x=x, y=y, col sep=comma] {dimacs-new/C2000-9/mssp-0.txt}; 
	\addplot[color=purple, line width=0.50mm] table [x=x, y=y, col sep=comma] {dimacs-new/C2000-9/sssp-0.txt}; 
	\addplot[color=blue, line width=0.50mm] table [x=x, y=y, col sep=comma] {dimacs-new/C2000-9/aco.txt}; 
	\addplot[color=gray, line width=0.50mm] table [x=x, y=y, col sep=comma] {dimacs-new/C2000-9/gurobi.txt}; 
	\addplot[color=black, line width=0.50mm, dashed] table [x=x, y=y, col sep=comma] {dimacs-new/C2000-9/lscc.txt};  
	\addplot[color=orange, line width=0.50mm, dashed] table [x=x, y=y, col sep=comma] {dimacs-new/C2000-9/fastwclq.txt};  
	\addplot[color=green, line width=0.50mm,dashed] table [x=x, y=y, col sep=comma] {dimacs-new/C2000-9/tsm.txt};
   \nextgroupplot[%
    title= \footnotesize C4000-5,
    ]
	\addplot[color=red, line width=0.50mm] table [x=x, y=y, col sep=comma] {dimacs-new/C4000-5/mssp-0.txt}; 
	\addplot[color=purple, line width=0.50mm] table [x=x, y=y, col sep=comma] {dimacs-new/C4000-5/sssp-0.txt}; 
	\addplot[color=blue, line width=0.50mm] table [x=x, y=y, col sep=comma] {dimacs-new/C4000-5/aco.txt}; 
	\addplot[color=gray, line width=0.50mm] table [x=x, y=y, col sep=comma] {dimacs-new/C4000-5/gurobi.txt}; 
	\addplot[color=black, line width=0.50mm, dashed] table [x=x, y=y, col sep=comma] {dimacs-new/C4000-5/lscc.txt};  
	\addplot[color=orange, line width=0.50mm, dashed] table [x=x, y=y, col sep=comma] {dimacs-new/C4000-5/fastwclq.txt};  
	\addplot[color=green, line width=0.50mm, dashed] table [x=x, y=y, col sep=comma] {dimacs-new/C4000-5/tsm.txt};
   \nextgroupplot[%
    title= \footnotesize hamming10-4,
    ]
	\addplot[color=red, line width=0.50mm] table [x=x, y=y, col sep=comma] {dimacs-new/hamming10-4/mssp-0.txt}; 
	\addplot[color=purple, line width=0.50mm] table [x=x, y=y, col sep=comma] {dimacs-new/hamming10-4/sssp-0.txt}; 
	\addplot[color=blue, line width=0.50mm] table [x=x, y=y, col sep=comma] {dimacs-new/hamming10-4/aco.txt}; 
	\addplot[color=gray, line width=0.50mm] table [x=x, y=y, col sep=comma] {dimacs-new/hamming10-4/gurobi.txt}; 
	\addplot[color=black, line width=0.50mm, dashed] table [x=x, y=y, col sep=comma] {dimacs-new/hamming10-4/lscc.txt};  
	\addplot[color=orange, line width=0.50mm, dashed] table [x=x, y=y, col sep=comma] {dimacs-new/hamming10-4/fastwclq.txt};  
	\addplot[color=green, line width=0.50mm, dashed] table [x=x, y=y, col sep=comma] {dimacs-new/hamming10-4/tsm.txt};
   \nextgroupplot[%
    title= \footnotesize keller6,
    ]
	\addplot[color=red, line width=0.50mm] table [x=x, y=y, col sep=comma] {dimacs-new/keller6/mssp-0.txt}; 
	\addplot[color=purple, line width=0.50mm] table [x=x, y=y, col sep=comma] {dimacs-new/keller6/sssp-0.txt}; 
	\addplot[color=blue, line width=0.50mm] table [x=x, y=y, col sep=comma] {dimacs-new/keller6/aco.txt};
	\addplot[color=gray, line width=0.50mm] table [x=x, y=y, col sep=comma] {dimacs-new/keller6/gurobi.txt}; 
	\addplot[color=black, line width=0.50mm, dashed] table [x=x, y=y, col sep=comma] {dimacs-new/keller6/lscc.txt};  
	\addplot[color=orange, line width=0.50mm, dashed] table [x=x, y=y, col sep=comma] {dimacs-new/keller6/fastwclq.txt};  
	\addplot[color=green, line width=0.50mm, dashed] table [x=x, y=y, col sep=comma] {dimacs-new/keller6/tsm.txt};
   \nextgroupplot[%
    title= \footnotesize MANN\_a45,
    ]
	\addplot[color=red, line width=0.50mm] table [x=x, y=y, col sep=comma] {dimacs-new/MANN-a45/mssp-0.txt}; 
	\addplot[color=purple, line width=0.50mm] table [x=x, y=y, col sep=comma] {dimacs-new/MANN-a45/sssp-0.txt}; 
	\addplot[color=blue, line width=0.50mm] table [x=x, y=y, col sep=comma] {dimacs-new/MANN-a45/aco.txt}; 
	\addplot[color=gray, line width=0.50mm] table [x=x, y=y, col sep=comma] {dimacs-new/MANN-a45/gurobi.txt}; 
	\addplot[color=black, line width=0.50mm, dashed] table [x=x, y=y, col sep=comma] {dimacs-new/MANN-a45/lscc.txt};  
	\addplot[color=orange, line width=0.50mm, dashed] table [x=x, y=y, col sep=comma] {dimacs-new/MANN-a45/fastwclq.txt};  
	\addplot[color=green, line width=0.50mm, dashed] table [x=x, y=y, col sep=comma] {dimacs-new/MANN-a45/tsm.txt};
   \nextgroupplot[%
    title= \footnotesize p\_hat1500-3,
    ]
	\addplot[color=red, line width=0.50mm] table [x=x, y=y, col sep=comma] {dimacs-new/p-hat1500-3/mssp-0.txt}; 
	\addplot[color=purple, line width=0.50mm] table [x=x, y=y, col sep=comma] {dimacs-new/p-hat1500-3/sssp-0.txt}; 
	\addplot[color=blue, line width=0.50mm] table [x=x, y=y, col sep=comma] {dimacs-new/p-hat1500-3/aco.txt}; 
	\addplot[color=gray, line width=0.50mm] table [x=x, y=y, col sep=comma] {dimacs-new/p-hat1500-3/gurobi.txt}; 
	\addplot[color=black, line width=0.50mm, dashed] table [x=x, y=y, col sep=comma] {dimacs-new/p-hat1500-3/lscc.txt};  
	\addplot[color=orange, line width=0.50mm, dashed] table [x=x, y=y, col sep=comma] {dimacs-new/p-hat1500-3/fastwclq.txt};  
	\addplot[color=green, line width=0.50mm, dashed] table [x=x, y=y, col sep=comma] {dimacs-new/p-hat1500-3/tsm.txt};
    \end{groupplot} 
    \node at (group c3r2.south) [anchor=south, yshift=-5cm, xshift=0.5cm] {\ref{grouplegend1}};       
	\end{tikzpicture}
\caption{The results generated by the generic and specialized methods on the MWCP using large graph instances in the DIMACS benchmark. The results of generic methods are in solid lines and those of specialized methods are in dashed lines. For a method, the mean statistic of best-found solutions is reported averaged over $10$ independent runs. We note that the ASP and SSSP are trained using a set of $26$ easy graphs with known optimal solutions.}
\label{fig:mwc}
\end{figure*}

Figure~\ref{fig:mwc} shows the results for the MWCP on $8$ large DIMACS benchmarks. Here, we focus on the comparison between several generic methods (i.e., those in solid lines). Overall, it can be seen that ASP significantly outperforms all other methods on all graphs, with an exception on MANN\_a45 where Gurobi is very competitive. Notably, ML-based methods (SSSP and ASP) can find (much) better solutions than ACO from the beginning in the solving process, and ACO may not find better solutions given longer computational time for some graphs such as `C4000-5' and `keller6'. When more computational time is given for certain graphs such as `hamming10-4', ACO can find better solutions than SSSP's due to its online-learning nature, but it falls much behind ASP using both online and offline learning. 

Figure~\ref{fig:mwc} also presents the results for three specialized methods or solvers (i.e., those in dashed lines) for the MWCP -- TSM~\citep{jiang2018two}, LSCC+BMS~\citep{wang2016two} and Fastwclq~\citep{cai2016fast}. We can observe that ASP is very competitive with these specialized methods in terms of solution quality, only slightly worse than LSCC+BMS. Although ASP does not achieve the best results, ASP can find many more near-optimal solutions efficiently than these specialized methods. This makes ASP advantageous to be used as a heuristic-pricing method for CG. As we will show in Section~\ref{sec:cgbp}, ASP can boost CG and B\&P for solving the graph coloring problem with the MWCP being the pricing problem.

\begin{figure*}[ht!]
    \centering
	\begin{tikzpicture}
	\begin{groupplot}[group style = {group size = 3 by 1, horizontal sep = 45pt, vertical sep = 45pt}, height=0.3\textwidth, width=0.3\textwidth, grid style={line width=.1pt, draw=gray!10},major grid style={line width=.2pt,draw=gray!30}, xlabel = \scriptsize Seconds, ylabel = \scriptsize Primal gap, y label style={at={(axis description cs:-0.1,.5)},anchor=south},x label style={at={(axis description cs:0.5,-0.05)},anchor=north}, xmajorgrids=true, ymajorgrids=true,  major tick length=0.05cm, minor tick length=0.0cm, legend style={column sep = 1pt,font=\scriptsize, legend columns = 4,draw=none}, title style={yshift=-1.5ex}]
   \nextgroupplot[%
    legend to name=grouplegend2,
    title=\footnotesize100 cities,
    ]
	\addplot[color=red, line width=0.30mm] table [x=x, y=y, col sep=comma] {figures/tsp/mssp_100.txt};\addlegendentry{ASP}
	\addplot[color=blue, line width=0.30mm] table [x=x, y=y, col sep=comma] {figures/tsp/aco_100.txt};\addlegendentry{ACS}
  	\addplot[color=purple, line width=0.30mm] table [x=x, y=y, col sep=comma] {figures/tsp/ortools_100.txt};\addlegendentry{OR-Tools}
	\addplot[color=black, line width=0.30mm] table [x=x, y=y, col sep=comma] {figures/tsp/gurobi_100.txt};\addlegendentry{Gurobi}
	\nextgroupplot[%
    title=\footnotesize300 cities,
    ]
	\addplot[color=red, line width=0.30mm] table [x=x, y=y, col sep=comma] {figures/tsp/mssp_300.txt};
	\addplot[color=blue, line width=0.30mm] table [x=x, y=y, col sep=comma] {figures/tsp/aco_300.txt};
         \addplot[color=purple, line width=0.30mm] table [x=x, y=y, col sep=comma] {figures/tsp/ortools_300.txt};
	\addplot[color=black, line width=0.30mm] table [x=x, y=y, col sep=comma] {figures/tsp/gurobi_300.txt};

	\nextgroupplot[%
    title=\footnotesize1000 cities,
    ]
	\addplot[color=red, line width=0.30mm] table [x=x, y=y, col sep=comma] {figures/tsp/mssp_1000.txt};
	\addplot[color=blue, line width=0.30mm] table [x=x, y=y, col sep=comma] {figures/tsp/aco_1000.txt};
      \addplot[color=purple, line width=0.30mm] table [x=x, y=y, col sep=comma] {figures/tsp/ortools_1000.txt};
	\addplot[color=black, line width=0.30mm] table [x=x, y=y, col sep=comma] {figures/tsp/gurobi_1000.txt};
    \end{groupplot} 
    \node at (group c1r1.north) [anchor=west, yshift=0.75cm, xshift=2.5cm] {\ref{grouplegend2}};           
	\end{tikzpicture}
\caption{The results on the TSP using ASP, ACS, OR-Tools, and Gurobi. For each problem size, the results are averaged over $100$ test instances.}
\label{fig:tsp-extended}
\end{figure*}
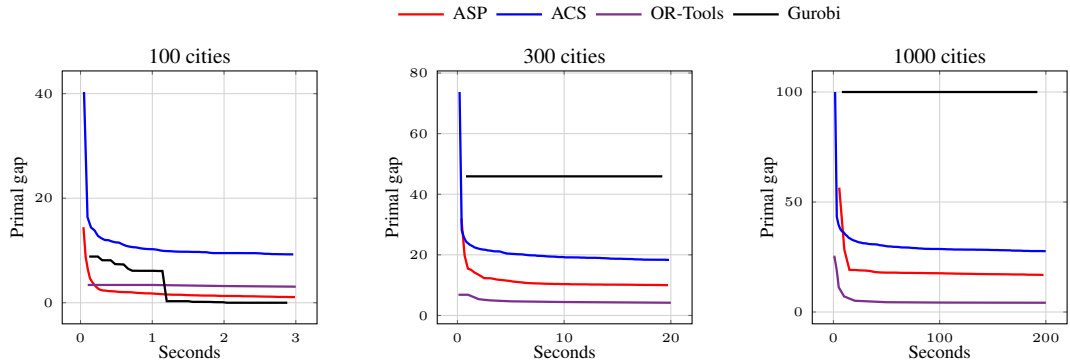

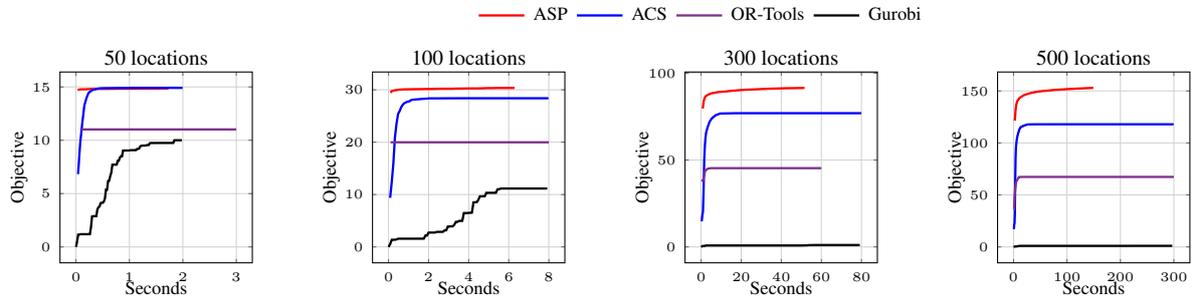
\begin{figure*}[ht!]
    \centering
	\begin{tikzpicture}
	\begin{groupplot}[group style = {group size = 4 by 1, horizontal sep = 45pt, vertical sep = 45pt}, height=0.25\textwidth, width=0.25\textwidth, grid style={line width=.1pt, draw=gray!10},major grid style={line width=.2pt,draw=gray!30}, xlabel = \scriptsize Seconds, ylabel = \scriptsize Objective, y label style={at={(axis description cs:-0.12,.5)},anchor=south},x label style={at={(axis description cs:0.5,-0.05)},anchor=north}, xmajorgrids=true, ymajorgrids=true,  major tick length=0.05cm, minor tick length=0.0cm, legend style={column sep = 1pt,font=\scriptsize, legend columns = 4,draw=none},  title style={yshift=-1.5ex}]
   \nextgroupplot[%
    ,legend to name=grouplegend,
    title=\footnotesize50 locations,
    ]
	\addplot[color=red, line width=0.30mm] table [x=x, y=y, col sep=comma] {figures/op/mssp_50unif.txt}; \addlegendentry{ASP}
	\addplot[color=blue, line width=0.30mm] table [x=x, y=y, col sep=comma] {figures/op/aco_50unif.txt}; \addlegendentry{ACS}
        \addplot[color=purple, line width=0.30mm] table [x=x, y=y, col sep=comma] {figures/op/ortools_unif50.txt};\addlegendentry{OR-Tools}
	\addplot[color=black, line width=0.30mm] table [x=x, y=y, col sep=comma] {figures/op/gurobi_50unif.txt}; \addlegendentry{Gurobi}
   \nextgroupplot[%
    title=\footnotesize100 locations,
    ]
	\addplot[color=red, line width=0.30mm] table [x=x, y=y, col sep=comma] {figures/op/mssp_100unif.txt};
	\addplot[color=blue, line width=0.30mm] table [x=x, y=y, col sep=comma] {figures/op/aco_100unif.txt};
        \addplot[color=purple, line width=0.30mm] table [x=x, y=y, col sep=comma] {figures/op/ortools_unif100.txt};
	\addplot[color=black, line width=0.30mm] table [x=x, y=y, col sep=comma] {figures/op/gurobi_100unif.txt};
	\nextgroupplot[%
    title=\footnotesize300 locations,
    ]
	\addplot[color=red, line width=0.30mm] table [x=x, y=y, col sep=comma] {figures/op/mssp_300unif.txt};
	\addplot[color=blue, line width=0.30mm] table [x=x, y=y, col sep=comma] {figures/op/aco_300unif.txt};
    \addplot[color=purple, line width=0.30mm] table [x=x, y=y, col sep=comma] {figures/op/ortools_unif300.txt};
	\addplot[color=black, line width=0.30mm] table [x=x, y=y, col sep=comma] {figures/op/gurobi_300unif.txt};
	\nextgroupplot[%
    title=\footnotesize500 locations,
    ]
	\addplot[color=red, line width=0.30mm] table [x=x, y=y, col sep=comma] {figures/op/mssp_500unif.txt};
	\addplot[color=blue, line width=0.30mm] table [x=x, y=y, col sep=comma] {figures/op/aco_500unif.txt};
        \addplot[color=purple, line width=0.30mm] table [x=x, y=y, col sep=comma] {figures/op/ortools_unif500.txt};
	\addplot[color=black, line width=0.30mm] table [x=x, y=y, col sep=comma] {figures/op/gurobi_500unif.txt};
    \end{groupplot} 

    \node at (group c1r1.north) [anchor=west, yshift=0.75cm, xshift=4cm] {\ref{grouplegend}};       
	\end{tikzpicture}
\caption{The results on the OP~(Maximization) using ASP, ACS, OR-Tools, and Gurobi. The test OP instances are of prize type `uniform'. For each problem size, the results are averaged over $100$ test instances. For ASP, the logistic model in ASP is optimized using OP instances with $500$ locations. }
\label{fig:op-extended}
\end{figure*}

Figures~\ref{fig:tsp-extended} and \ref{fig:op-extended} show the results for the two routing problems, the TSP and OP. As compared with ACO, it can be seen that ASP can find much better solutions very efficiently, consistent with the observations on the MWCP. The other generic method, Gurobi, is only effective in tackling small instances of these problems. For these routing problems, we have also tested a heuristic routing solver in OR-Tools, which is an open-source software suite for optimization\footnote{OR-Tools: \url{https://developers.google.com/optimization}}. The heuristic solver is configured using recommended settings. Specifically, a constructive heuristic is used to generate an initial feasible solution, starting from a route ``start" node and then connecting it to the node which produces the cheapest route segment. This solution is then iteratively refined by an improvement heuristic, named guided local search~\citep{voudouris2010guided}. The results show that OR-Tools achieves better performance than ASP on the TSP (Figure~\ref{fig:tsp-extended}) but obtains much worse results than ASP on the OP (Figure~\ref{fig:op-extended}). Lastly, we note that highly specialized solvers exist for these well-studied problems, such as LKH-3 for the TSP and Compass~\citep{kobeaga2018efficient} for the OP. Such a solver incorporates a number of specialized methods, and ASP is not expected to be competitive with these methods.

\subsection{ASP using different ML models}
\label{subsec:exp-other}

In this part, we show that our proposed ASP framework is generic -- it can improve the prediction quality for several classes of ML models including a classification tree and a feed-forward neural network (FFNN). The case study is conducted on the MWCP using the DIMACS benchmarks. The classification tree is trained using extreme gradient boosting~\citep{friedman2001greedy, chen2015xgboost}, referred to as XGBT. The FFNN is implemented using the Keras package~\citep{keras2022} and configured as follows. We set the number of layers to $3$, the number of neurons for the hidden layer to $32$, the activation function for the input layer and the hidden layer to Rectified Linear Units~\citep{nair2010rect}, the learning algorithm to Adam~\citep{kingma2015adam} with its default parameters, and the loss function to the cross-entropy loss function~\citep{bishop2006pattern}.

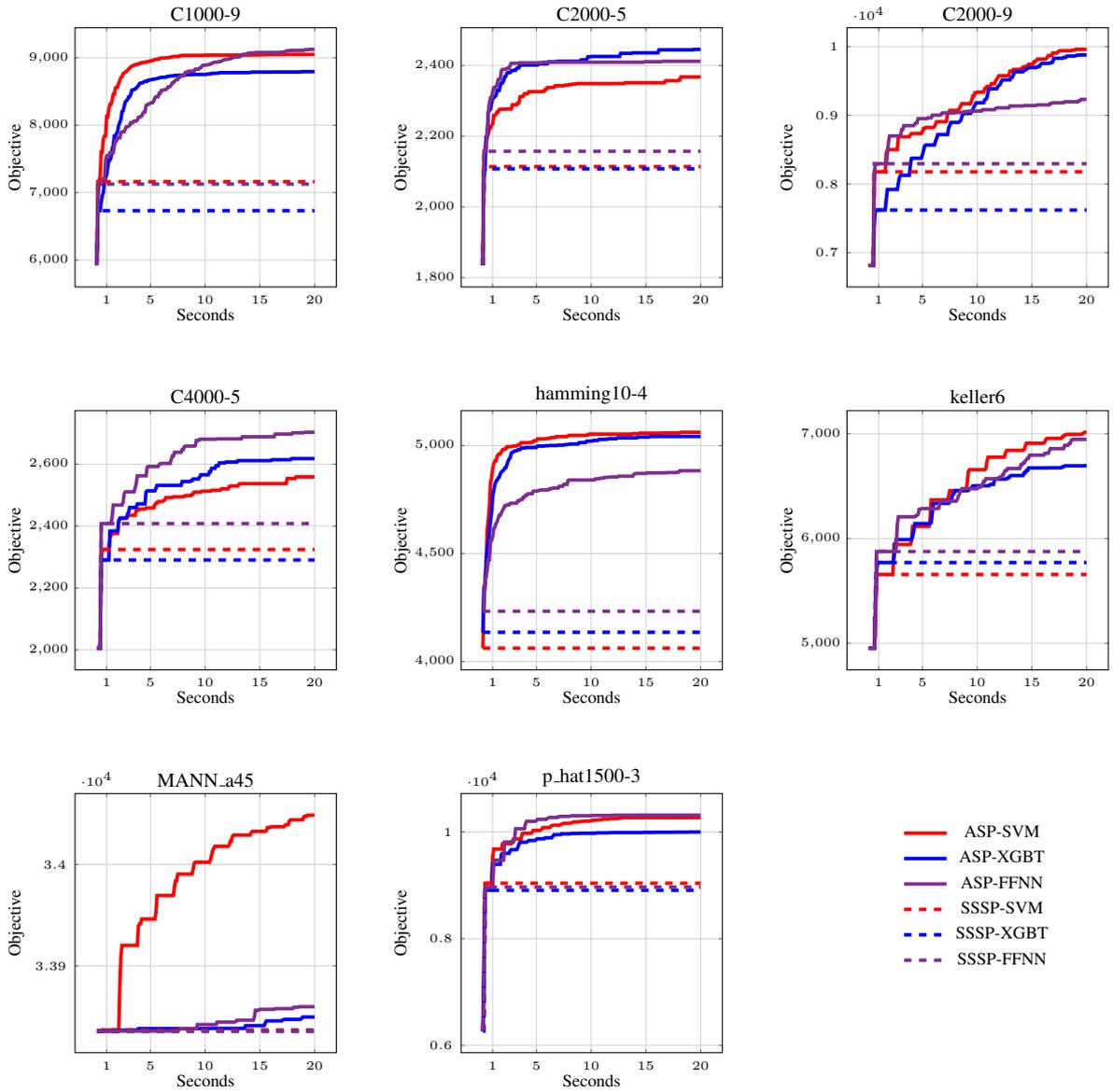
\begin{figure}[ht!]
    \centering
	\begin{tikzpicture}
	\begin{groupplot}[group style = {group size = 3 by 3, horizontal sep = 50pt, vertical sep = 50pt}, height=0.32\textwidth, width=0.32\textwidth, grid style={line width=.1pt, draw=gray!10},major grid style={line width=.2pt,draw=gray!30}, xlabel = \scriptsize Seconds, ylabel = \scriptsize Objective, y label style={at={(axis description cs:-0.16,.5)},anchor=south}, xtick = {1,5,10,15,20},x label style={at={(axis description cs:0.5,-0.05)},anchor=north}, xmajorgrids=true, ymajorgrids=true,  major tick length=0.05cm, minor tick length=0.0cm, legend style={column sep = 1pt,font=\scriptsize, legend columns = 1,draw=none},title style={yshift=-1.5ex}]
  \nextgroupplot[%
  legend to name=grouplegend3,
    title= \footnotesize C1000-9,
    ]

	\addplot[color=red, line width=0.50mm] table [x=x, y=y, col sep=comma] {dimacs-new/C1000-9/mssp-0.txt}; \addlegendentry{ASP-SVM}
	\addplot[color=blue, line width=0.50mm] table [x=x, y=y, col sep=comma] {dimacs-new/C1000-9/mssp-2.txt}; \addlegendentry{ASP-XGBT}
	\addplot[color=purple, line width=0.50mm] table [x=x, y=y, col sep=comma] {dimacs-new/C1000-9/mssp-3.txt}; \addlegendentry{ASP-FFNN}
	\addplot[color=red, line width=0.50mm, dashed] table [x=x, y=y, col sep=comma] {dimacs-new/C1000-9/sssp-0.txt}; \addlegendentry{SSSP-SVM}
	\addplot[color=blue, line width=0.50mm, dashed] table [x=x, y=y, col sep=comma] {dimacs-new/C1000-9/sssp-2.txt}; \addlegendentry{SSSP-XGBT}
	\addplot[color=purple, line width=0.50mm, dashed] table [x=x, y=y, col sep=comma] {dimacs-new/C1000-9/sssp-3.txt}; \addlegendentry{SSSP-FFNN}
  \nextgroupplot[%
    title= \footnotesize C2000-5,
    ]
	\addplot[color=red, line width=0.50mm] table [x=x, y=y, col sep=comma] {dimacs-new/C2000-5/mssp-0.txt}; 
	\addplot[color=blue, line width=0.50mm] table [x=x, y=y, col sep=comma] {dimacs-new/C2000-5/mssp-2.txt}; 
	\addplot[color=purple, line width=0.50mm] table [x=x, y=y, col sep=comma] {dimacs-new/C2000-5/mssp-3.txt}; 
	\addplot[color=red, line width=0.50mm, dashed] table [x=x, y=y, col sep=comma] {dimacs-new/C2000-5/sssp-0.txt}; 
	\addplot[color=blue, line width=0.50mm, dashed] table [x=x, y=y, col sep=comma] {dimacs-new/C2000-5/sssp-2.txt}; 
	\addplot[color=purple, line width=0.50mm, dashed] table [x=x, y=y, col sep=comma] {dimacs-new/C2000-5/sssp-3.txt}; 
  \nextgroupplot[%
    title= \footnotesize C2000-9,
    ]
	\addplot[color=red, line width=0.50mm] table [x=x, y=y, col sep=comma] {dimacs-new/C2000-9/mssp-0.txt}; 
	\addplot[color=blue, line width=0.50mm] table [x=x, y=y, col sep=comma] {dimacs-new/C2000-9/mssp-2.txt}; 
	\addplot[color=purple, line width=0.50mm] table [x=x, y=y, col sep=comma] {dimacs-new/C2000-9/mssp-3.txt}; 
	\addplot[color=red, line width=0.50mm, dashed] table [x=x, y=y, col sep=comma] {dimacs-new/C2000-9/sssp-0.txt}; 
	\addplot[color=blue, line width=0.50mm, dashed] table [x=x, y=y, col sep=comma] {dimacs-new/C2000-9/sssp-2.txt}; 
	\addplot[color=purple, line width=0.50mm, dashed] table [x=x, y=y, col sep=comma] {dimacs-new/C2000-9/sssp-3.txt}; 
  \nextgroupplot[%
    title= \footnotesize C4000-5,
    ]
	\addplot[color=red, line width=0.50mm] table [x=x, y=y, col sep=comma] {dimacs-new/C4000-5/mssp-0.txt}; 
	\addplot[color=blue, line width=0.50mm] table [x=x, y=y, col sep=comma] {dimacs-new/C4000-5/mssp-2.txt}; 
	\addplot[color=purple, line width=0.50mm] table [x=x, y=y, col sep=comma] {dimacs-new/C4000-5/mssp-3.txt}; 
	\addplot[color=red, line width=0.50mm, dashed] table [x=x, y=y, col sep=comma] {dimacs-new/C4000-5/sssp-0.txt}; 
	\addplot[color=blue, line width=0.50mm, dashed] table [x=x, y=y, col sep=comma] {dimacs-new/C4000-5/sssp-2.txt}; 
	\addplot[color=purple, line width=0.50mm, dashed] table [x=x, y=y, col sep=comma] {dimacs-new/C4000-5/sssp-3.txt}; 
  \nextgroupplot[%
    title= \footnotesize hamming10-4,
    ]
	\addplot[color=red, line width=0.50mm] table [x=x, y=y, col sep=comma] {dimacs-new/hamming10-4/mssp-0.txt}; 
	\addplot[color=blue, line width=0.50mm] table [x=x, y=y, col sep=comma] {dimacs-new/hamming10-4/mssp-2.txt}; 
	\addplot[color=purple, line width=0.50mm] table [x=x, y=y, col sep=comma] {dimacs-new/hamming10-4/mssp-3.txt};
	\addplot[color=red, line width=0.50mm, dashed] table [x=x, y=y, col sep=comma] {dimacs-new/hamming10-4/sssp-0.txt}; 
	\addplot[color=blue, line width=0.50mm, dashed] table [x=x, y=y, col sep=comma] {dimacs-new/hamming10-4/sssp-2.txt}; 
	\addplot[color=purple, line width=0.50mm, dashed] table [x=x, y=y, col sep=comma] {dimacs-new/hamming10-4/sssp-3.txt};
  \nextgroupplot[%
    title= \footnotesize keller6,
    ]
	\addplot[color=red, line width=0.50mm] table [x=x, y=y, col sep=comma] {dimacs-new/keller6/mssp-0.txt}; 
	\addplot[color=blue, line width=0.50mm] table [x=x, y=y, col sep=comma] {dimacs-new/keller6/mssp-2.txt}; 
	\addplot[color=purple, line width=0.50mm] table [x=x, y=y, col sep=comma] {dimacs-new/keller6/mssp-3.txt};
	\addplot[color=red, line width=0.50mm, dashed] table [x=x, y=y, col sep=comma] {dimacs-new/keller6/sssp-0.txt}; 
	\addplot[color=blue, line width=0.50mm, dashed] table [x=x, y=y, col sep=comma] {dimacs-new/keller6/sssp-2.txt}; 
	\addplot[color=purple, line width=0.50mm, dashed] table [x=x, y=y, col sep=comma] {dimacs-new/keller6/sssp-3.txt};
  \nextgroupplot[%
    title= \footnotesize MANN\_a45,
    ]
	\addplot[color=red, line width=0.50mm] table [x=x, y=y, col sep=comma] {dimacs-new/MANN-a45/mssp-0.txt}; 
	\addplot[color=blue, line width=0.50mm] table [x=x, y=y, col sep=comma] {dimacs-new/MANN-a45/mssp-2.txt}; 
	\addplot[color=purple, line width=0.50mm] table [x=x, y=y, col sep=comma] {dimacs-new/MANN-a45/mssp-3.txt};
	\addplot[color=red, line width=0.50mm, dashed] table [x=x, y=y, col sep=comma] {dimacs-new/MANN-a45/sssp-0.txt}; 
	\addplot[color=blue, line width=0.50mm, dashed] table [x=x, y=y, col sep=comma] {dimacs-new/MANN-a45/sssp-2.txt}; 
	\addplot[color=purple, line width=0.50mm, dashed] table [x=x, y=y, col sep=comma] {dimacs-new/MANN-a45/sssp-3.txt};
  \nextgroupplot[%
    title= \footnotesize p\_hat1500-3,
    ]

	\addplot[color=red, line width=0.50mm] table [x=x, y=y, col sep=comma] {dimacs-new/p-hat1500-3/mssp-0.txt}; 
	\addplot[color=blue, line width=0.50mm] table [x=x, y=y, col sep=comma] {dimacs-new/p-hat1500-3/mssp-2.txt}; 
	\addplot[color=purple, line width=0.50mm] table [x=x, y=y, col sep=comma] {dimacs-new/p-hat1500-3/mssp-3.txt};
	\addplot[color=red, line width=0.50mm, dashed] table [x=x, y=y, col sep=comma] {dimacs-new/p-hat1500-3/sssp-0.txt}; 
	\addplot[color=blue, line width=0.50mm, dashed] table [x=x, y=y, col sep=comma] {dimacs-new/p-hat1500-3/sssp-2.txt}; 
	\addplot[color=purple, line width=0.50mm, dashed] table [x=x, y=y, col sep=comma] {dimacs-new/p-hat1500-3/sssp-3.txt};
    \end{groupplot} 
    \node at (group c3r2.south) [anchor=south, yshift=-4.5cm, xshift=0cm] {\ref{grouplegend3}};       
	\end{tikzpicture}
\caption{On the MWCP, the results for ASP and SSSP with different ML models for solution prediction.}
\label{fig:mwc2}
\end{figure}

Figure~\ref{fig:mwc2} shows the results for ASP and SSSP using different ML models, respectively. Here, the SSSP refers to the special case of ASP that predicts the optimal solution only once before the search. We highlight the following two observations.  Firstly, for the same ML model, ASP can improve over SSSP by a large margin, confirming the efficacy and the general applicability of ASP for different ML models. Secondly, for the same approach, i.e., either ASP or SSSP, there is no clear winner among the different ML models used within that approach. We have chosen the linear-SVM as our default ML model for ASP, because of its very fast prediction times. This can be crucial when ASP is repeatedly used for solving pricing problems within branch-and-price.

%% file: 05-cgbp.tex
\section{ASP for Column Generation and Branch-and-Price}
\label{sec:cgbp}

	


This section uses ASP as a pricing heuristic for CG~\citep{lubbecke2005selected}, to boost the exact B\&P method~\citep{barnhart1998branch} for solving the graph coloring problem (GCP). CG aims to capture the columns in the optimal solution of a linear program~(LP), and these columns are gradually generated by solving a sequence of NP-hard pricing problems. Using ASP as the pricing heuristic for solving the pricing problems can accelerate CG in capturing optimal LP columns, because ASP can generate many high-quality columns efficiently. This section is built on our previous work~\citep{shen2022enhancing}, which proposes a SSSP approach as the pricing heuristic for CG and the exact B\&P method for solving GCP. 

\subsection{Background}



In GCP, the goal is to assign a minimum number of colors to vertices in a graph, such that every pair of the adjacent vertices does not share the same color~\citep{malaguti2010survey}. To solve GCP, B\&P recursively decomposes the original GCP (root node) into subproblems (children nodes) and prunes a node safely without further expansion if its lower bound was no better than the current best-found solution. For GCP, a tighter lower bound can be obtained by solving the LP relaxation of its set-covering formulation (GCP-SC), defined as using a minimum number of maximal independent sets~(MISs) to cover all the vertices in a graph~\citep{mehrotra1996column}. Let the binary variable~$z_\mathcal{I}$ indicate whether a MIS~$\mathcal{I} \in \mathbb{I}$ is used to cover the vertices in a graph $G(\mathcal{V},\mathcal{E})$, where $\mathbb{I}$ is the set of all the possible MISs; $\mathcal{V}$ denotes the set of vertices; $\mathcal{E}$ denotes the set of edges. GCP-SC can be formulated as
\begin{align}
    \min_{\bm{z}} \;& \sum_{m \in \mathcal{I}} z_\mathcal{I}, & \text{(GCP-SC)}\\
        s.t. \;& \sum_{\mathcal{I} \in \mathbb{I}, i\in \mathcal{I}} z_\mathcal{I} \geq 1; &  i \in \mathcal{V}, \label{eq:at-least-one} \\
               & z_\mathcal{I} \in \{0,1\}; &\mathcal{I} \in \mathbb{I},
\end{align}
\noindent where Constraints (24) specify that, for each vertex $i$, at least one MIS containing that vertex $i \in \mathcal{I}$ must be in use.

Given GCP-SC, its LP relaxation can be obtained by relaxing integer constraints, i.e., replacing Constraints (25) to $z_\mathcal{I} \in [0,1]; \mathcal{I} \in \mathbb{I}$. Solving the LP can yield a tight lower bound that can help B\&P aggressively prune nodes and potentially reduce its computational time. However, solving such a LP itself can be challenging to existing methods (e.g., the simplex methods or the interior point methods~\citep{dantzig2016linear}), because the LP may need an exponential number of variables (or columns) to represent all MISs in a graph.

CG is an iterative method for solving large-scale LPs, based on the fact that an optimal LP solution typically has only a tiny proportion of optimal LP columns, i.e., columns with non-zero values. CG starts by solving a restricted master problem (RMP) that considers a subset of columns (or MISs) in the original LP and gradually includes new MISs that can improve the solution to the current RMP, i.e., those with negative reduced costs (NRC). Searching for NRC MISs can be formulated as an optimization problem (i.e., the pricing problem), with minimizing the reduced cost as the objective function, defined for GCP-SC as
\begin{align}
    \min_{\bm{x}} \;& 1 - \sum_{i \in \mathcal{V}} \pi_i\cdot x_i, & \text{(MWISP)}\\
        s.t. \;& x_i + x_j \leq 1; & (i,j)\in \mathcal{E}, \\
              & x_i \in \{0,1\}; & i \in \mathcal{V},
\end{align}
\noindent where $x_i$ denotes whether a vertex $i$ is a part of the MIS and $\pi_i$ holds the optimal dual solution associated with that vertex $i$ to the current RMP. As can be seen, the pricing problem is the maximum weight independent problem (MWISP) on the graph $G(\mathcal{V},\mathcal{E}, \mathcal{W})$ with the optimal dual solution as vertex weights $\mathcal{W} \coloneqq \{\pi_i\;|\;i \in \mathcal{V}\}$. Note that the MWISP is equivalent to the MWCP on its complementary graph $G(\mathcal{V},\overline{\mathcal{E}}, \mathcal{W})$ as shown in Equations~(19)-(21).

CG alternates between solving the RMP and the MWISP, and the computational process terminates when no NRC column exists. This certifies that the RMP has already captured all the optimal LP columns to the original LP; hence the optimal solution to the RMP is also optimal to the original LP.

\subsection{Experiments for Column Generation} \label{subsubsec:ret-cg}
We test CG using different pricing heuristics, for solving the LP relaxation of GCP-SC. Specifically, at an iteration of CG, we use a heuristic method to tackle the NP-hard MWISP aiming to quickly find NRC MISs, and add to the RMP at most $n$ NRC MISs in increasing order of their reduced costs. When the heuristic method fails to find any NRC MIS, an exact method, TSM~\citep{jiang2018two}, is used to optimally solve the pricing problem. If the optimal MIS has a negative reduced cost, it is included in the RMP to start the next iteration. Note that the setup of a pricing heuristic and the column-selection strategy are common practices for CG~\citep{lubbecke2005selected}.

We use the ASP implementation for the MWCP described in Section~\ref{subsec:mwcp}. To collect the training data, we use CG to solve LPs on a set of $10$ graphs and record the optimally solved MWISPs encountered during the execution of CG. The graphs are from the graph coloring benchmarks\footnote{\url{https://sites.google.com/site/graphcoloring/files}}, including: ``3-FullIns\_4", ``queen12\_12", ``1-Insertions\_6", ``mug88\_25", ``DSJC125.5", ``flat300\_20\_0", ``flat300\_26\_0", ``DSJC1000.9", ``DSJC250.1", and ``queen11\_11". These graphs are selected under two considerations. Firstly, their pricing problems are easy to solve, i.e., collecting the training data is easy. Secondly, they are of different characteristics so that a ML model can learn common knowledge that is applicable to different unseen graphs. 

For the compared pricing methods, we include SSSP~\citep{shen2022enhancing}, ACO~\citep{xu2007improved}, Gurobi~\citep{gurobi2018gurobi} and specialized methods to the MWCP. These specialized methods include TSM~\citep{jiang2018two}, LSCC+BMS~\citep{wang2016two}, and Fastwclq~\citep{cai2016fast}. For each method, we set a maximum execution time of $30$ seconds for solving a MWISP instance, to prevent one from spending too much time at a particular CG iteration. Note that the exact methods, Gurobi and TSM, are also subject to the cutoff time and treated as `heuristics'. Furthermore, we classify the compared methods into several classes, and the termination criterion for each class is tuned empirically to maximize the performance of CG~\citep{shen2022enhancing}. In particular, the sampling-based methods (i.e., ASP, SSSP, ACO, and Fastwclq) construct a set of $50n$ MISs at a CG iteration, where $n$ is the number of vertices in a graph. We set ASP to perform $50$ iterations with $n$ MISs generated at each iteration. 

Our test problem instances are generated using the graph coloring benchmarks. A problem instance is a RMP initialized with $10n$ randomly generated columns (or MISs). We generate $24$ problem instances with different random seeds $[1,\cdots,24]$ on a set of $81$ test graphs, resulting in a total number of $1944$ test problem instances. These graphs are used for instance generation, because their initial RMPs can be solved within a reasonable time. That is $500$ seconds under an overall cutoff time of $1800$ seconds in our case. Note that we use the default LP solver in Gurobi~\citep{gurobi2018gurobi} for solving RMPs in CG. All of our experiments are conducted on a server with $32$ CPUs (AMD EPYC Processor, 2245 MHz) and $128$ GB RAM.

\begin{figure}[ht!]
	\begin{tikzpicture}
	\begin{groupplot}[group style = {group size = 2 by 1, horizontal sep = 60pt}, height=0.475\textwidth, width=0.475\textwidth, grid style={line width=.1pt, draw=gray!10},major grid style={line width=.2pt,draw=gray!30}, xlabel = \small Time in Seconds, xtick = {0,600,1200,1800}, ticklabel style = {font=\small}, xmajorgrids=true, ymajorgrids=true,  major tick length=0.05cm, minor tick length=0.0cm, legend style={font=\small, column sep = 1pt, legend columns = 4,draw=none}]
   \nextgroupplot[%
    legend to name=group1,
    ylabel = \small \# Solve Instances,
    y label style={at={(axis description cs:-0.12,.5)},anchor=south},
    ]

	\addplot[color=red, line width=0.45mm] table [x=x, y=y, col sep=comma] {cg/small/solving-curve/mssp-cs0.txt};\addlegendentry{\small CG-ASP}
	\addplot[color=blue, line width=0.45mm] table [x=x, y=y, col sep=comma] {cg/small/solving-curve/sssp-cs0.txt};\addlegendentry{\small CG-SSSP}
	\addplot[color=purple, line width=0.45mm] table [x=x, y=y, col sep=comma] {cg/small/solving-curve/aco.txt};\addlegendentry{\small CG-ACO}
	\addplot[color=orange, line width=0.45mm] table [x=x, y=y, col sep=comma] {cg/small/solving-curve/fastwclq.txt};\addlegendentry{\small CG-Fastwclq}
	\addplot[color=gray, line width=0.45mm] table [x=x, y=y, col sep=comma] {cg/small/solving-curve/lscc.txt};\addlegendentry{\small CG-LSCC}
	\addplot[color=olivegreen, line width=0.45mm] table [x=x, y=y, col sep=comma] {cg/small/solving-curve/tsm.txt};\addlegendentry{\small CG-TSM}
	\addplot[color=black, line width=0.45mm] table [x=x, y=y, col sep=comma] {cg/small/solving-curve/gurobi.txt};\addlegendentry{\small CG-Gurobi}
   \nextgroupplot[%
    ylabel = \small LP Obective,
    y label style={at={(axis description cs:-0.12,.5)},anchor=south},
    ]
    

	\addplot[color=red, line width=0.45mm] table [x=x, y=y, col sep=comma] {cg/small/lp-curve/mssp-cs0.txt};
	\addplot[color=blue, line width=0.45mm] table [x=x, y=y, col sep=comma] {cg/small/lp-curve/sssp-cs0.txt};
	\addplot[color=purple, line width=0.45mm] table [x=x, y=y, col sep=comma] {cg/small/lp-curve/aco.txt};
	\addplot[color=orange, line width=0.45mm] table [x=x, y=y, col sep=comma] {cg/small/lp-curve/fastwclq.txt};
	\addplot[color=gray, line width=0.45mm] table [x=x, y=y, col sep=comma] {cg/small/lp-curve/lscc.txt};
	\addplot[color=olivegreen, line width=0.45mm] table [x=x, y=y, col sep=comma] {cg/small/lp-curve/tsm.txt};
	\addplot[color=black, line width=0.45mm] table [x=x, y=y, col sep=comma] {cg/small/lp-curve/gurobi.txt};
    \end{groupplot} 
    \node at (group c1r1.north) [anchor=north, yshift=1.5cm, xshift=4cm] {\pgfplotslegendfromname{group1}}; 
	\end{tikzpicture}
\caption{Results for CG with different pricing methods. \textbf{Left:} the number of solved instances out of $1944$ test instances. \textbf{Right:} the objective values of the RMP (the lower the better), averaged over all problem instances using the geometric mean.}
\label{fig:ret_small}
\end{figure}
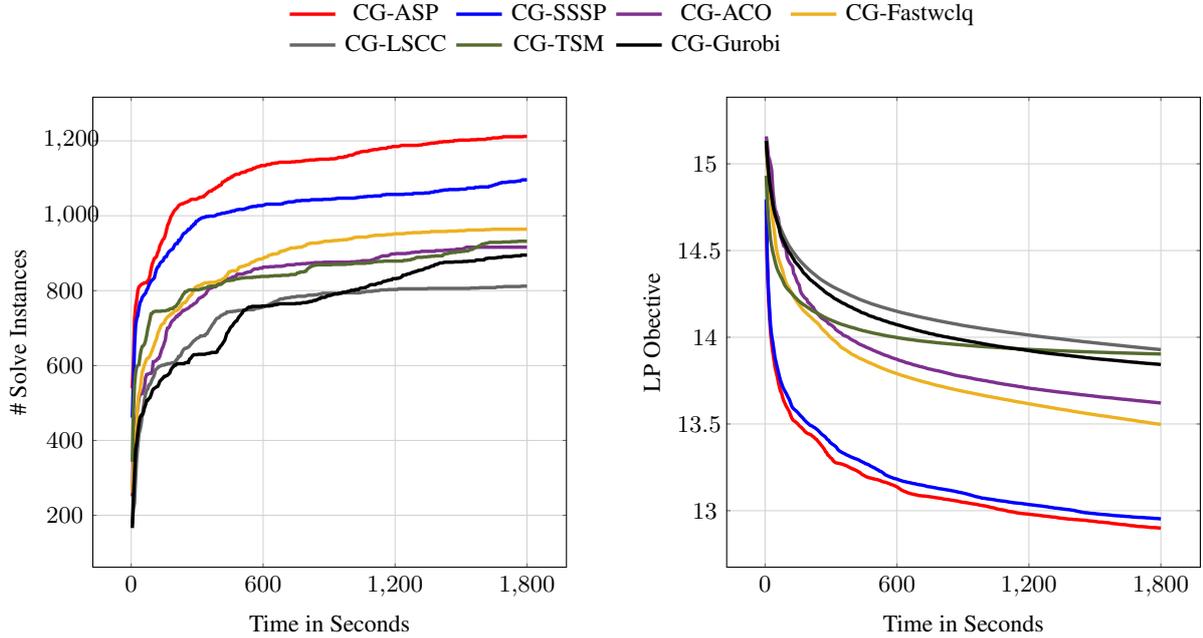

Figure~\ref{fig:ret_small} shows the results for CG using different pricing heuristics. In the left sub-figure, we report the total number of solved problem instances within $1800$ seconds. In the right sub-figure, we report the progress of CG with respect to the LP objective value, which takes the unsolved problem instances into account as well. It is interesting to observe the following: 1) CG with the ML-based pricing methods (ASP or SSSP) can make more substantial progress than CG with other methods. This is because ML-based pricing heuristics can generate \emph{many high-quality} columns efficiently. Both the number of NRC columns and the quality of the columns are crucial to the computational efficiency of CG, as noted by previous studies~\citep{lubbecke2005selected}. In at least one of the aspects, ML-based heuristic methods are superior to other pricing methods~\citep{shen2022enhancing}. For example, the performance of CG with LSCC is not competitive because LSCC cannot find a large number of NRC MISs, though LSCC is shown to be very efficient in finding a single (near-)optimal solution~(Figure~\ref{fig:mwc}). 2) A comparison of the ML-based heuristic methods shows that ASP further improves the performance of CG compared to SSSP. Specifically, CG with ASP solves many more problem instances than CG with SSSP at any given runtime, as shown in the left sub-figure. This is because ASP can find better-quality MISs than SSSP for solving MWISPs, as evidenced by the results in Figure~\ref{fig:mwc}. However, the results in the right sub-figure of Figure~\ref{fig:ret_small} show that CG-ASP is marginally better than CG-SSSP in terms of the mean LP objective values. These observations reveal an interesting phenomenon in CG called the tailing-off effect, i.e., the difficulty of CG in converging to the optimal LP solution given near-optimal ones. 

\subsection{Experiments for Branch-and-Price}

We test B\&P using CG with different pricing heuristics, for solving the integer problem GCP. We use the B\&P code from an academic mixed-integer-programming solver, SCIP~\citep{gamrath2020scip}. In particular, for solving the pricing problems, the B\&P code uses an efficient greedy search~\citep{mehrotra1996column} as the pricing heuristic and an exact method called $t$-clique. Besides that, this B\&P implementation incorporates some specialized techniques for solving GCPs from previous studies~\citep{mehrotra1996column,malaguti2011exact}, such as a branching method that operates on the original compact formulation of GCP, a local search method that finds the initial primal solution of GCP (i.e., an upper bound) to be used as initial LP columns for RMP, an early branching method that may terminate CG early and safely by comparing Lagrangian lower bound of the current RMP with its objective value. In addition, the B\&P can also benefit from the functionalities provided by the SCIP solver, such as primal heuristics for finding better feasible solutions and column-management algorithms for efficiently re-optimizing RMPs during the CG process. These components form a reasonable B\&P implementation for testing different pricing heuristics. Moreover, we note that additional techniques may further improve the performance of B\&P, such as primal heuristics specialized to CG~\citep{JoncourMSSV10} or methods for stabilizing dual solutions of RMPs~\citep{du1999stabilized}.  

We refer to the default setup of B\&P as B\&P-def, and compare it with the B\&P variants that replace the greedy search to one of our ML-based heuristic pricing methods (i.e., ASP or SSSP), respectively referred to as B\&P-ASP and B\&P-SSSP. B\&P-SSSP uses a sample size of $10n$ by default~\citep{shen2022enhancing}. To ensure a fair comparison between the ML-based pricing methods, we use B\&P-ASP to construct the same number of samples evenly divided into $10$ ASP iterations. From newly generated NRC MISs, we add at most $\theta$ MISs into the RMP in the increasing order of their reduced costs, $\theta=n$ for the root node and $\theta=0.1n$ for children nodes. This is because typically a child node already contains a large number of useful columns inherent from its parent. Hence including a large number of columns is rarely helpful but increases the computational overhead of solving subsequent RMPs.

For each method, we perform $1584$ runs on a set of~$66$ graphs in the graph coloring benchmarks. The excluded graphs are either too easy (all methods can solve them within $10$ seconds) or too hard to solve (all methods cannot solve the LP at the root node) under the cutoff time of $8000$ seconds. When the LP at the root node is solved to optimality, we report the duality gap, defined as 
\begin{equation}
\text{duality gap} \coloneqq 100\% \times \frac{upper\_bound - global\_lower\_bound}{upper\_bound},   
\end{equation}

\noindent where the upper bound is the objective value of the best-found solution and the global lower bound is the smallest lower bound amongst the remaining open tree nodes. 

\begin{table}[ht!]
    \caption{Results for B\&P with different pricing methods. For each graph, a method solves $24$ problem instances (generated with different random columns in the initial root RMP) and the average results are reported. The solving time is an average over optimally solved GCPs, and `N/A' is placed if a method does not solve any instances on that graph in 8000 seconds. Similarly, the duality gap for a graph is an average over instances whose root LP has been solved to optimality, and `N/A' is placed if the root LP cannot be solved in any run.}
    \label{tab:bp}
    \centering
    \resizebox{\linewidth}{!}{\begin{tabular}{@{}clrr|rrr|rrr|rrr|rrr@{}}
        \toprule
       \multirow{2}{*}{Group} & \multirow{2}{*}{Graph} & \multirow{2}{*}{\# nodes} & \multirow{2}{*}{Density} & \multicolumn{3}{c}{\# optimally Solved} & \multicolumn{3}{c}{Solving time (sec)} & \multicolumn{3}{c}{\# root LP solved} & \multicolumn{3}{c}{Duality gap (\%)} \\
    & & & & ASP & SSSP & Greedy & ASP & SSSP & Greedy & ASP & SSSP & Greedy & ASP & SSSP & Greedy\\
    \cmidrule(lr){1-1}\cmidrule{2-4}\cmidrule(lr){5-7}\cmidrule(lr){8-10}\cmidrule(lr){11-13}\cmidrule(lr){14-16}
      
    \multirow{13}{*}{1}  & DSJR500.5 & 486 & 0.972 & \textbf{24} & 17 & 15 & 786.3 & 1370.9 & 782.4 & 24 & 24 & 24 & \textbf{0.0} & 0.2 & 0.3 \\
     & 1-FullIns\_4 & 38 & 0.364 & \textbf{24} & 18 & 23 & \textbf{411.7} & 1289.2 & 675.0  & 24 & 24 & 24 & \textbf{0.0} & 5.0 & 0.8 \\
     & queen9\_9 & 81 & 0.652 & \textbf{24} & 0 & 2 & \textbf{12.6} & N/A & 2204.6 & 24 & 24 & 24 & \textbf{0.0} & 10.0 & 9.2 \\
     & DSJC125.5 & 125 & 0.502 & \textbf{2} & 0 & 0 & \textbf{4905.3} & N/A & N/A & 24 & 24 & 24 & \textbf{9.9} & 15.4 & 15.8 \\
     &  2-Insertions\_3 & 37 & 0.216 & \textbf{24} & 17 & 9 & 4105.3 & 4571.5 & 3448.5 & 24 & 24 & 24 & \textbf{0.0} & 7.3 & 15.6 \\  
     & ash608GPIA & 1215 & 0.021 & 24 & 24 & 0 & \textbf{2840.8} & 4899.6 & N/A & 24 & 24 & 0 & 0.0 & 0.0 & N/A \\
     & will199GPIA & 660 & 0.054 & 24 & 24 & 24 & \textbf{54.8} & 70.4 & 128.9 & 24 & 24 & 24 & 0.0 & 0.0 & 0.0 \\
     & 1-Insertions\_4 & 67 & 0.210 & 0 & 0 & 0 & N/A & N/A & N/A & 24 & 24 & 24 & \textbf{20.0} & 37.5 & 40.0 \\
     & myciel5 & 47 & 0.437 & 0 & 0 & 0 & N/A & N/A & N/A & 24 & 24 & 24 & \textbf{16.7} & 22.9 & 32.6 \\
     & DSJC250.9 & 250 & 0.896 & 0 & 0 & 0 & N/A & N/A & N/A & 24 & 24 & 24 & \textbf{2.1} & 2.3 & 4.1 \\
     & r1000.1c & 709 & 0.969 & 0 & 0 & 0 & N/A & N/A & N/A & 24 & 24 & 24 & \textbf{5.9} & 6.2 & 7.1 \\
     & 1-FullIns\_5 & 78 & 0.277 & 0 & 0 & 0 & N/A & N/A & N/A & 24 & 24 & 24 & \textbf{16.7} & 31.9 & 33.3 \\
     & myciel6 & 95 & 0.338 & 0 & 0 & 0 & N/A & N/A & N/A & 24 & 24 & 24 & \textbf{28.6} & 31.0 & 42.9 \\ 
    
    \midrule
    \multirow{6}{*}{2}  & le450\_5c & 450 & 0.194 & 13 & \textbf{21} & 1 & 4263.3 & 4024.2 & 2009.9 & 13 & \textbf{21} & 1 & 0.0 & 0.0 & 0.0 \\
    & le450\_5d & 450 & 0.193 & 17 & \textbf{21} & 10 & 4483.0 & 2905.7 & 2126.4 & 21 & \textbf{23} & 11 & 3.2 & 1.5 & 1.5 \\
    & flat300\_26\_0 & 300 & 0.965 & 24 & 24 & 24 & 361.6 & \textbf{223.0} & 1910.6 & 24 & 24 & 24 & 0.0 & 0.0 & 0.0 \\  
    & school1\_nsh & 326 & 0.547 & 24 & 24 & 24 & 685.8 & \textbf{627.0} & 759.5 & 24 & 24 & 24 & 0.0 & 0.0 & 0.0 \\
    & ash331GPIA & 661 & 0.038 & 24 & 24 & 24 & 59.1 & \textbf{53.1} & 307.4 & 24 & 24 & 24 & 0.0 & 0.0 & 0.0 \\
    & queen16\_16 & 256 & 0.387 &  0 & 0  & 0  & N/A & N/A & N/A & 23 & \textbf{24} & 0 & 11.1 & 11.1 & N/A \\

    \midrule
    \multirow{10}{*}{3}  & flat300\_20\_0 & 300 & 0.953 & 24 & 24 & 24 & 84.9 & 84.9 & 343.6 & 24 & 24 & 24 & 0.0 & 0.0 & 0.0 \\
    & school1 & 355 & 0.603 & 24 & 24 & 24 & 51.7 & 52.0 & 1966.9 & 24 & 24 & 24 & 0.0 & 0.0 & 0.0 \\
    & DSJC1000.9 & 1000 & 0.900 &  0 & 0  & 0  & N/A & N/A & N/A & 24 & 24 & 17 & 11.5 & 11.5 & 11.5 \\
    & le450\_25d & 433 & 0.366 &  0 & 0  & 0  & N/A & N/A & N/A & 24 & 24 & 0 & 10.7 & 10.7 & N/A \\ 
    & le450\_25c & 435 & 0.362 &  0 & 0  & 0  & N/A & N/A & N/A & 24 & 24 & 0 & 10.7 & 10.7 & N/A \\
    & le450\_15b & 410 & 0.187 &  0 & 0  & 0  & N/A & N/A & N/A & 24 & 24 & 0 & 6.3 & 6.3 & N/A \\
    &  qg.order40 & 1600 & 0.098 &  0 & 0  & 0  & N/A & N/A & N/A & 24 & 24 & 0 & 2.4 & 2.4 & N/A \\
    & le450\_15a & 407 & 0.189 &  0 & 0  & 0  & N/A & N/A & N/A & 24 & 24 & 0 & 6.3 & 6.3 & N/A \\
    & wap06a & 703 & 0.288 &  0 & 0  & 0  & N/A & N/A & N/A & 24 & 24 & 0 & 4.8 & 4.8 & N/A \\
    & queen15\_15 & 225 & 0.411 &  0 & 0  & 0  & N/A & N/A & N/A & 24 & 24 & 0 & 11.8 & 11.8 & N/A \\

    \midrule
    \multirow{5}{*}{4} & 3-Insertions\_5 & 1406 & 0.020 &  0 & 0  & 0  & N/A & N/A & N/A & 1 & 0 & \textbf{24} & 50.0 & N/A & 50.0 \\
    & 2-Insertions\_5 & 597 & 0.044 &  0 & 0  & 0  & N/A & N/A & N/A & 1 & 1 & \textbf{24} & 50.0 & 50.0 & 50.0 \\
    & 4-Insertions\_4 & 475 & 0.032 &  0 & 0  & 0  & N/A & N/A & N/A & 23 & 15 & \textbf{24} & 40.0 & 40.0 & 40.0 \\    
    & r1000.5 & 966 & 0.477 &  0 & 0  & 0  & N/A & N/A & N/A & 24 & 24 & 24 & 9.3 & 7.8 & \textbf{1.3} \\
    & 1-Insertions\_5 & 202 & 0.121 &  0 & 0  & 0  & N/A & N/A & N/A & 24 & 24 & 24 & 34.0 & 43.8 & \textbf{33.3} \\
        \bottomrule
    \end{tabular}}
\end{table}

Table~\ref{tab:bp} shows the results for B\&P. According to the performances of B\&P using different pricing heuristics, the test graphs are divided into four groups and are described respectively. Group~1 includes $13$ graphs where B\&P with ASP outperforms others, and the improvements are substantial on several graphs. Specifically, B\&P-ASP solves GCP on `queen9\_9' in all runs with an average time of $12.6$ seconds and on `DSJC125.5' in two runs, whereas other methods cannot solve any with a cutoff time of $8000$ seconds; B\&P-ASP reduces the duality gap almost by half compared to other methods for the graphs `1-Insertions\_4' and `1-FullIns-5’. Notably, B\&P-ASP also significantly improves on B\&P-SSSP. This is expected as B\&P repeatedly invokes CG at each node of the B\&P process, and CG-ASP is more efficient than CG-SSSP in general (see Section~\ref{subsubsec:ret-cg}). Group 2 consists of $6$ graphs on which B\&P with SSSP outperforms others. In particular, B\&P-SSSP is considerably better than B\&P-ASP on the two graphs with the prefix `le450-5'. For these graphs, the ML predictions by ASP may have converged fast to local optima, resulting in a less diverse set of sampled MISs. Group~3 shows the results on $10$ graphs, where the performances of B\&P with the two ML-based pricing heuristics are comparable and are much better than B\&P with its default greedy search. This is because for these graphs CG can be sped up substantially by including multiple MISs at each iteration. Both ML-based pricing methods can generate many high-quality MISs, whereas the default greedy heuristic can only generate a single MIS per CG iteration. In contrast, the results in Group~4 show that B\&P with the ML-based pricing methods do not outperform B\&P-def. This is because, for these graphs, the exact method can efficiently solve the pricing problems already, and competitive performance of CG can be obtained by disabling the use of a heuristic method~\citep{shen2022enhancing}. Lastly, we note that on some test graphs the results for the three methods are comparable, and further experimental results can be found in the Supplementary Materials.

%% file: 06-conclusion.tex
\section{Conclusion}
\label{sec:conclusion}

This paper proposes an adaptive solution prediction (ASP) framework that iteratively refines feature representations of decision variables using statistical information extracted from best-found solutions in the search process (i.e., online learning), enabling an offline-trained ML model to predict the optimal solutions in an adaptive manner. As the heuristic search (i.e., probabilistic sampling) continuously supplies better-quality solutions, more accurate statistical information can be obtained to better characterize these variables. As a result, an offline-trained ML model can better predict the optimal solution, which in turn better guides the search.

The efficacy of ASP is demonstrated empirically. Firstly, ASP is generic: 1) it can be used to improve the prediction quality of different ML models to find better solutions, as evidenced by a case study on the maximum weight clique problem (MWCP); 2) ASP can be potentially adapted to solve different COPs, and this study considers MWCP and two routing problems. Secondly, ASP is competitive: 1) compared with conventional generic heuristic methods (without ML), ASP can find better solutions much more efficiently because of its offline learning via ML; 2) compared with existing ML-based methods, ASP can find much better solutions especially for large problem instances, benefiting from its online learning; 3) compared with highly-specialized methods on MWCP, ASP can find many more high-quality solutions efficiently. This makes it advantageous to be used as a heuristic-pricing method to boost branch-and-price (B\&P) for solving the graph coloring problem. 

In future work, we would like to test B\&P with ASP for heuristic pricing on more combinatorial optimization problems such as vehicle routing problems. Our goal is to develop a \emph{generic} ML-based pricing heuristic to speed up CG and B\&P for solving the Dantzig-Wolfe reformulation of combinatorial optimization problems.